\documentclass[pmlr]{jmlr}
\makeatletter
\let\Ginclude@graphics\@org@Ginclude@graphics
\def\ps@jmlrtps{\let\@mkboth\@gobbletwo
  \def\@oddhead{}\let\@evenhead\@oddhead
  \def\@oddfoot{\hfill\small\rmfamily\thepage\hfill}\let\@evenfoot\@oddfoot}
\makeatother
\newcommand{\coltauthor}[1]{\author{#1}}

\title[Polylogarithmic Regret in Degenerate LQR]{Polylogarithmic Regret via Certainty Equivalence in Degenerate LQR}

\usepackage{times}

\coltauthor{
 \Name{Luca Vignola} \Email{lvignola@ethz.ch}\\
 \Name{Melanie Zeilinger} \Email{mzeilinger@ethz.ch}\\
 \Name{Bruce D. Lee} \Email{bruce.lee@ethz.ch}\\
 \addr ETH Zürich, Zürich, Switzerland%
}

\usepackage{amsmath,amssymb,amsfonts}
\usepackage{graphicx}
\usepackage{epstopdf}
\usepackage{algorithm}
\usepackage{textcomp}
\usepackage{mathtools}
\usepackage{etoolbox}
\usepackage[T1]{fontenc}
\usepackage[utf8]{inputenc}
\usepackage{dsfont}

\usepackage[english]{babel}
\usepackage{amsfonts}
\usepackage{bbm}
\usepackage{nicefrac}
\usepackage{xfrac}
\usepackage{booktabs}
\usepackage{multirow}
\usepackage{tikz}
\usepackage{tikz-3dplot}
\usetikzlibrary{3d,arrows.meta,patterns,decorations.markings,calc,positioning}
\usepackage{wrapfig}
\usepackage{microtype}      %
\usepackage{needspace}
\usepackage{placeins}
\usepackage{bm}

\usepackage{empheq}
\usepackage{color}
\definecolor{cool-blue}{HTML}{1D6996}
\definecolor{cool-purple}{HTML}{5F4690}
\definecolor{cool-teal}{HTML}{38A6A5}
\definecolor{cool-green}{HTML}{0F8554}
\definecolor{dark-green}{rgb}{0.0,0.5,0.0}

\usepackage[capitalize]{cleveref}
\let\For\relax
\let\If\relax
\usepackage[noend]{algpseudocode}

\AtBeginEnvironment{lemma}{\crefalias{theorem}{lemma}}
\AtBeginEnvironment{corollary}{\crefalias{theorem}{corollary}}
\AtBeginEnvironment{proposition}{\crefalias{theorem}{proposition}}
\AtBeginEnvironment{definition}{\crefalias{theorem}{definition}}
\AtBeginEnvironment{remark}{\crefalias{theorem}{remark}}
\renewenvironment{proof}[1][\proofname]
  {\par\noindent{\bfseries\upshape #1\ }}
  {\jmlrQED}

\crefname{definition}{Definition}{Definitions}
\newtheorem{assumption}{Assumption}
\crefname{assumption}{Assumption}{Assumptions}
\crefname{theorem}{Theorem}{Theorems}
\crefname{corollary}{Corollary}{Corollaries}
\crefname{lemma}{Lemma}{Lemmas}
\crefname{remark}{Remark}{Remarks}

\usepackage[theorems,skins]{tcolorbox}
\usepackage{stackengine}

\stackMath          %

\definecolor{bruce-maroon}{RGB}{128,0,0}

\usepackage{todonotes}

\begin{document}

\maketitle

\begin{abstract}%
Singular input penalties, as in minimum variance control, can reduce what must be identified to control an unknown linear system.
In online linear quadratic control, certainty equivalence can achieve polylogarithmic regret when the optimal controller is identifiable from data collected under optimal feedback.
With a positive definite input penalty, however, this condition is restrictive, since for an invertible state matrix it requires all unknown parameters to be identifiable under optimal feedback.
We show that singular input penalties relax this requirement, expanding the class of problems in which certainty equivalence achieves polylogarithmic regret.
Specifically, we prove that a certainty-equivalence algorithm with a short initial exploration phase achieves $O(\log^2 T)$ expected regret over $T$ interactions whenever controller identifiability holds.
We complement this guarantee with an $\Omega(\sqrt T)$ minimax lower bound for a family containing a system where controller identifiability fails.
For minimum variance control, knowing the range of the input matrix suffices for controller identifiability, a weaker requirement than knowing the full input matrix, as assumed in prior work with positive definite input penalties.
This condition also covers underactuated systems, as we illustrate with a cart--pole example.
\end{abstract}


\section{Introduction}
Learning to control an unknown system during operation requires balancing exploration, which collects informative data to reduce uncertainty about the dynamics, with exploitation, which uses the available knowledge to achieve good control performance.
The linear--quadratic regulator (LQR) with unknown dynamics has become a standard benchmark for studying this trade-off in online control and in reinforcement learning with continuous state and action spaces.
In LQR, under a fixed feedback $u_t=Kx_t$, the next state $x_{t+1}=Ax_t+Bu_t+w_t$  depends on $(A,B)$ only through the closed loop $A+BK$, 
so models that share this closed loop cannot be distinguished from the collected data.
A \emph{certainty-equivalence} (CE) learner, which always plays the controller that is optimal for a nominal system estimate, may therefore keep collecting data that cannot reveal such model errors. %
Moreover, classical results show that convergence of self-tuning LQ controllers to the optimal controller generally requires identification of the true parameter~\citep{Polderman86}. This motivated a variety of exploration strategies~\citep{pmlr-v19-abbasi-yadkori11a,mania2019extended,bartosOptimisticOnlineLQR2026}, which attain the worst-case optimal $\widetilde{\mathcal O}(\sqrt{T})$ regret rate over $T$ interactions with the unknown system~\citep{pmlr-v119-simchowitz20a}.
CE can achieve faster, polylogarithmic rates when additional structure, such as prior knowledge of the dynamics, makes the data generated under optimal feedback sufficiently informative to identify all unknown parameters. 
In this work, we show that singular input penalties (\emph{degenerate} LQR) can allow CE to achieve finite-time polylogarithmic expected regret even when the dynamics are not fully identifiable under optimal feedback.

\subsection{Related Work}
\paragraph{Exploration in online LQR.}
In the finite-time regret formulation of adaptive LQR, \citet{pmlr-v19-abbasi-yadkori11a} obtained the first
$\widetilde{\mathcal O}(\sqrt{T})$ regret bound using optimism in the face of
uncertainty, and \citet{cohenLearningLinearQuadraticRegulators2019} subsequently gave the first
computationally efficient algorithm attaining the same rate.
\citet{mania2019extended} established this rate for CE with injected input noise, %
and \citet{pmlr-v119-simchowitz20a} proved that the regret rate achieved by such naive exploration cannot be improved in the worst case for fully unknown systems, up to logarithmic factors.
These approaches provide $\widetilde{\mathcal O}(\sqrt{T})$ regret guarantees, but do not capture the faster rates that may be attainable when the data generated under optimal feedback are sufficiently informative.

\paragraph{CE under structured dynamics.}
With suitable structure in the dynamics, the data generated by a CE controller can be sufficiently informative for learning. CE can then achieve improved regret rates without persistent exploration.
For example, a greedy CE algorithm achieves polylogarithmic expected regret
when the input matrix $B_\star$ is known, or when the state matrix $A_\star$ is known
and the optimal feedback gain is sufficiently non-degenerate
\citep{casselLogarithmicRegretLearning2020,jedra2021minimalexpectedregretlinear}. \citet{faradonbeh2020adaptive} showed that an episodic CE scheme attains $O(\log^2T)$ asymptotic regret almost surely, assuming that the parameters lie in a known set over which the closed-loop prediction error uniformly bounds the error of the optimal gain.
\citet{lee2024nonasymptoticregretanalysisadaptive} established finite-time polylogarithmic expected-regret bounds for CE with a known linear representation of the dynamics when the optimal controller persistently excites the unknown parameters.

When some parameter directions cannot be detected from data collected under optimal feedback, the consequences depend on whether they affect the optimal controller. For known $A_\star$, \citet{ziemann2021uninformative} established a local asymptotic minimax lower bound of order $\sqrt T$ around uninformative instances, where optimal feedback fails to reveal a parameter direction that changes the optimal gain.
Motivated by this result, \citet{liu2026controller} showed
 that it suffices to identify the optimal gain rather than the dynamics, and used \emph{controller identifiability} (CI),
  which excludes such directions, to separate logarithmic from square-root local minimax regret.
These results assume a positive definite input penalty $R\succ0$. In this setting, CI requires all unknown parameters to be identifiable under optimal feedback whenever $A_\star$ is nonsingular~\citep[Corollary~2]{liu2026controller}. The extent to which singular input penalties relax this requirement remains unaddressed.

\paragraph{CE for degenerate cost functions.}
\looseness=-1 A separate line of work studies minimum variance control (MVC). %
In a related multivariate stochastic regression setting, \citet{LAI86} obtained asymptotically logarithmic regret by combining CE with sparse probing inputs.
\citet{laiWei1987SelfTuning} extended this approach to MVC of finite-order scalar autoregressive models with exogenous inputs (ARX models).
More recently, \citet{pmlr-v235-singh24b} established finite-time polylogarithmic regret for self-tuning regulation of ARX models, of order $\log T$ for bounded noise and $\log^3T$ for sub-Gaussian noise.
\citet{LIN198577} further argued that CE can fail to converge to the optimal controller for cost criteria other than minimum output variance. Together, these results suggest
that degenerate input costs can also contribute to making full parameter identification unnecessary.
However, the MVC guarantees for ARX models concern scalar-output regulation without an input penalty. They do not provide guarantees for general quadratic state costs in underactuated systems, or for costs that penalize only some input directions.

\subsection{Contributions}
We extend controller identifiability to LQR with potentially singular input penalties, $R\succeq0$, and establish finite-time polylogarithmic expected-regret guarantees for CE. Specifically:
\begingroup
\csappto{@listi}{%
    \setlength{\topsep}{3pt}%
    \setlength{\itemsep}{2pt}%
    \setlength{\parsep}{0pt}%
}
\begin{itemize}
    \item In \Cref{sec:ci}, we restate CI in terms of the curvature of the cost, which makes the role of the cost function explicit, and extend it to $R\succeq0$. With $m$ inputs, a singular $R$ admits up to $m\dim\ker R$ additional parameter directions that cannot be detected from data collected under optimal feedback but are harmless for control, so that CI no longer requires full parameter identifiability, even for nonsingular $A_\star$.
    \item In \Cref{section:mvc_upperbound}, we prove that CE with a short warm-up achieves $O(\log^2T)$ expected regret whenever CI holds at the true system. The bound is non-asymptotic and requires only a stabilizing controller and coarse bounds on the system parameters. %
	\item In \Cref{section:mvc_lowerbound}, we prove a local $\Omega(\sqrt T)$ minimax lower bound around a system where CI fails, for any input penalty, including $R=0$.
    \item In \Cref{sec:mvc}, we specialize these results to MVC and show that CI holds whenever $\operatorname{range}(B_\star)$ is known. This prior is weaker than knowing $B_\star$ and covers underactuated systems, as we illustrate on a cart--pole example.
\end{itemize}
\endgroup

\section{Problem Setting}\label{sec:pb_setting}
We use $\|\cdot\|$ for the Euclidean vector or spectral matrix norm,
$\|\cdot\|_F$ for the Frobenius norm, and $\langle X,Y\rangle_F=\operatorname{Tr}(X^\top Y)$.
\subsection{Degenerate LQR}
\paragraph{Dynamics and model class.} We consider the linear system
\begin{equation}
    x_{t+1}=A_{\star} x_t + B_{\star} u_t + w_t,\qquad t\ge1,\qquad x_1=0,
\end{equation}
with i.i.d.\ process noise $w_t\sim\mathcal{N}(0,\sigma^2I_n)$ of known variance $\sigma^2>0$. The matrices $\Theta_\star=[A_\star\ B_\star]$, with $A_\star\in\mathbb{R}^{n\times n}$ and $B_\star\in\mathbb{R}^{n\times m}$, are unknown, but the learner is given an affine model class that contains them,
\begin{equation}\label{eq:a1_affine_model}
\begin{gathered}
 \mathcal C=\{\Theta(\theta):\theta\in\mathbb R^p\},\qquad
 \Theta(\theta)=\Theta_{\mathrm{ref}}+\mathcal E\theta,\\
 \mathcal E\theta=\textstyle\sum_{j=1}^{p}\theta_j E_j,\qquad
 \langle E_i,E_j\rangle_F=\mathbf 1_{\{i=j\}},
\end{gathered}
\end{equation}
where $\Theta_{\mathrm{ref}}$ and $E_j$ are known and $\Theta_\star=\Theta(\theta_\star)$ for an unknown
 $\theta_\star\in\mathbb R^p$. This class includes fully
  unknown dynamics, but can also encode prior structural knowledge. Any finite-dimensional affine family admits this representation after removing
    redundant coordinates and orthonormalizing its directions. In this representation, $\|\mathcal Eh\|_F=\|h\|_2$ and
 we write $[\Delta A_h\ \Delta B_h]=\mathcal Eh$ for $h\in\mathbb R^p$.

 \paragraph{Cost and optimal controller.} Policies are causal; that is, $u_t$ is chosen based on $x_1,\dots,x_t$, past inputs,
and internal randomization that is independent of the process noise.
Given known $Q\succ 0$ and $R\succeq 0$, the stage cost at step $t$ is
\begin{equation}
    c_t=x_t^\top Qx_t+u_t^\top Ru_t.
\end{equation}
For a model $\Theta=[A\ B]$ representing the linear dynamics and a linear controller $K$,
 we then define the following infinite-horizon average cost 
\begin{align}\label{eq:cost_definition}
    J_{\Theta}(K)
    \coloneqq
    \limsup\limits_{T\to\infty}\frac1T
    \mathbb E^{\Theta}_K\!\left[\sum_{t=1}^{T}c_t\right],
\end{align}
where the expectation is over the noise $w_t$ under the dynamics $\Theta$ and the feedback $u_t=Kx_t$. We omit this dependence of the expectation when clear from the context and write $J(K)\coloneqq J_{\Theta_\star}(K)$.

The following assumption ensures that the LQR problem and optimal controller are well defined even when $R$ is singular.
\begin{assumption}\label{assumption:riccati}
$(A_{\star},B_{\star})$ is stabilizable, $Q\succ0$, $R\succeq0$, and $\ker R\cap\ker B_\star=\{0\}$.
\end{assumption}

\looseness=-1 \Cref{thm:riccati_solution} shows that, under \Cref{assumption:riccati}, the discrete-time algebraic Riccati equation (DARE)
\[
 P=Q+A_\star^\top P A_\star
 -A_\star^\top P B_\star
 (R+B_\star^\top P B_\star)^{-1}
 B_\star^\top P A_\star
\]
admits a unique positive definite stabilizing solution $P_\star\succ0$. The same proposition proves that $u_t=K_\star x_t$ is the unique stationary linear feedback minimizing $J(K)$, with
\[
 S_\star=R+B_\star^\top P_\star B_\star\succ0,\qquad
 K_\star=-S_\star^{-1}B_\star^\top P_\star A_\star,
\]
and ensures that the closed-loop matrix $L_\star=A_\star+B_\star K_\star$ is Schur stable.
The assumption $Q\succ0$ ensures that every nonzero state direction is penalized, and the kernel condition that every unpenalized input direction affects the state; for $R=0$, the latter means that $B_\star$ has full column rank.

\paragraph{Regret and certainty equivalence.} Let $J_\star\coloneqq J_{\Theta_\star}(K_\star)$ denote the optimal infinite-horizon average cost. We measure the cost of learning over $T$ decision rounds by the cumulative regret
\begin{equation}
    R_T=\sum_{t=1}^T\left(c_t-J_\star\right),
\end{equation}
where $c_t$ denotes the cost incurred by the learner on the true system $\Theta_\star$. Our performance criterion is then %
the expected cumulative regret $\mathbb{E}[R_T]$, where the expectation is taken over
the process noise and any internal randomization of the learning algorithm.
\begin{remark}[Expected regret]
    Since $J_\star$ is an average over the noise, whereas each $c_t$ depends on its realization, the standard deviation of $R_T$ is of order $\sqrt T$, even under the optimal controller. For this definition of regret, polylogarithmic rates are therefore only meaningful in expectation. Asymptotic almost-sure polylogarithmic rates, such as those of~\citet{faradonbeh2020adaptive}, are instead stated for a regret that compares the cost of the learner with that of the optimal controller under the same noise realization, so that this contribution of the noise cancels. We also note that, compared to an analysis on a high-probability event, bounding the expected regret additionally requires controlling the cost of rare failures, such as estimates that destabilize the system, which can be very large~\citep{casselLogarithmicRegretLearning2020}.
\end{remark}

Finally, for an estimated model $\hat\Theta=[\hat A\ \hat B]$, we define the certainty-equivalent controller by
\begin{equation}\label{controller:CE}
    K(\hat\Theta)=\arg\min_{K} J_{\hat\Theta}(K).
\end{equation}
The minimizer is unique whenever $(\hat A,\hat B)$ satisfies the conditions of \Cref{assumption:riccati}
 in place of $(A_\star,B_\star)$ (\Cref{thm:riccati_solution}). For models in $\mathcal C$, we overload notation by writing $K(\theta)\coloneqq K(\Theta(\theta))$.

\subsection{Controller Identifiability}\label{sec:ci}
Fix a feedback gain $K$ that stabilizes the model $\Theta(\theta)=[A(\theta)\ B(\theta)]$. Under this feedback, let $p_\theta^K(y\mid x)$ denote the conditional density of the next state $y$ given the current state $x$. This density is Gaussian with mean $(A(\theta)+B(\theta)K)x$ and covariance $\sigma^2I_n$.
We quantify the information about $\theta$ provided by each transition using the stationary per-step Fisher information:
\[
 \mathcal I_\theta(K)=\mathbb E_{\theta,K}\!\left[
 \nabla_\theta\log p_\theta^K(y\mid x)
 \bigl(\nabla_\theta\log p_\theta^K(y\mid x)\bigr)^\top\right],
\]
where the expectation is over two consecutive states $(x,y)$ in the stationary closed loop and the gradient is taken with respect to $\theta$, holding $K$ fixed.
We further define
\[
 \mathcal I_\star=\mathcal I_{\theta_\star}(K_\star),\qquad
 \mathcal H_\star=
 \left.\nabla_\theta^2 J_{\Theta_\star}(K(\theta))
 \right|_{\theta=\theta_\star},
\]
where $\mathcal H_\star$ measures the local sensitivity of the CE controller's performance to identification errors. For $\theta$ sufficiently close to $\theta_\star$, $K(\theta)$ stabilizes the true system $\Theta_\star$, and, since $K(\theta_\star)=K_\star$ is optimal, a second-order Taylor expansion gives
\[
 J_{\Theta_\star}(K(\theta))-J_\star
 =\tfrac12(\theta-\theta_\star)^\top\mathcal H_\star(\theta-\theta_\star)
 +o(\|\theta-\theta_\star\|^2).
\]
We refer to the parameter directions in $\ker\mathcal I_\star$ as \emph{invisible} directions, because perturbations along them do not change the distribution of data collected under the fixed controller $K_\star$. We call directions in $\ker\mathcal H_\star$ \emph{harmless}, meaning that identification errors along them have no second-order effect on the CE controller's cost. Controller identifiability (CI) requires every invisible direction to be harmless.
\begin{definition}
\label{def:controller_identifiability}
Controller identifiability holds at $\theta_\star$ if
$\ker\mathcal I_\star\subseteq\ker\mathcal H_\star$.
\end{definition}

CI was introduced by \citet{liu2026controller}, who defined it for $R\succ0$ 
   by requiring the derivative of the optimal gain to vanish along invisible directions, thereby excluding the \emph{uninformative} instances of \citet{ziemann2021uninformative}. We show in \Cref{lemma:controller_identifiability_tests} that our definition agrees with that of \citet{liu2026controller} for $R\succ0$, while extending CI to potentially singular input penalties $R\succeq0$. %
	Our formulation directly highlights the role of the cost function in determining whether an invisible direction is harmless for performance.
	This leads to the following generalization of \citet[Corollary 2]{liu2026controller}, which quantifies how many additional harmless invisible directions a singular input penalty allows.

\begin{lemma}
\label{lemma:degenerate_cost_dimension}
Suppose \Cref{assumption:riccati} holds and CI holds at
$\theta_\star$. In the nonredundant affine representation
\eqref{eq:a1_affine_model},
\[
 \dim\ker\mathcal I_\star
 \leq m\bigl(n-\operatorname{rank}A_\star+\dim\ker R\bigr).
\]
\end{lemma}

The proof is given in \Cref{sec:ci_appendix}. When $A_\star$ is nonsingular and $R\succ0$, 
 CI therefore implies $\mathcal I_\star\succ0$, i.e., all unknown parameters must be
  identifiable under optimal feedback. A singular $R$ instead allows CI with up to $m\dim\ker R$ 
  invisible directions. %
  In \Cref{sec:mvc}, we give an MVC model class ($R=0$) for which CI holds with $\dim\ker\mathcal I_\star=m^2=m\dim\ker R$, matching the upper bound when $A_\star$ is nonsingular.
  
  \paragraph{Example.} Consider the scalar system $x_{t+1}=ax_t+bu_t+w_t$, where $a$ and $b$ are unknown and nonzero. Under optimal feedback $u_t=K_\star x_t$, the observations identify only the closed-loop coefficient $a+bK_\star$. The invisible directions therefore satisfy $\Delta a+K_\star\Delta b=0$. For $R>0$, \Cref{lemma:degenerate_cost_dimension} implies that CI would require full parameter identifiability because $a\ne0$, so CI fails. When $R=0$, $K_\star=-a/b$, and small invisible perturbations scale $a$ and $b$ by the same factor. These perturbations preserve the ratio $a/b$, hence the optimal gain and the CE controller's cost, so they are harmless. Thus CI holds even though $a$ and $b$ cannot be separately identified under optimal feedback.

\section{Main Results}
In this section, we show that, when CI holds at $\theta_\star$, a CE algorithm achieves polylogarithmic regret after a short exploration phase. \Cref{section:mvc_lowerbound} complements this result with a
worst-case lower bound over a family of systems including a plant where CI fails.
\subsection{Regret Upper Bound} \label{section:mvc_upperbound}

\subsubsection{Algorithm}\label{sec:prior_information}
Similarly to \citet{casselLogarithmicRegretLearning2020}, the learner is supplied with a gain $K_0$ that stabilizes the true system. It first plays $K_0$ with
independent exploration noise $\eta_t\sim\mathcal N(0,\sigma^2I_m)$
for $\tau_0-1$ rounds, where $\tau_0$ is an input to the algorithm (\Cref{alg:mvc}). The learner then proceeds in epochs, each of length at most $\tau_0$.
At each epoch boundary, the learner computes a regularized least-squares estimate restricted to the model class $\mathcal C$
\begin{equation}\label{eq:constrained_ridge}
 \hat\Theta_t=\operatorname*{arg\,min}_{\Theta\in\mathcal C}
 \left\{\sum_{s=1}^{t-1}\|x_{s+1}-\Theta z_s\|_2^2
                         +\lambda\|\Theta\|_F^2\right\},
\end{equation}
where $z_t=(x_t^\top,u_t^\top)^\top$ and the regularization parameter $\lambda$ is an input to the algorithm. The learner plays the CE gain $K(\hat\Theta_t)$ until the next epoch boundary, without additional exploration noise. The algorithm permanently falls back to $K_0$ if the estimated
model has no positive definite stabilizing Riccati solution with
$\hat S_t=R+\hat B_t^\top\hat P_t\hat B_t\succ0$, or if $\|K(\hat\Theta_t)\|>\kappa$ or $\|x_t\|^2>x_b$, where the thresholds $\kappa, x_b$ are also inputs to the algorithm. These checks prevent a poor estimate from destabilizing the closed loop, which is needed to bound the regret in expectation. Our analysis uses coarse bounds on the problem parameters to construct explicit values of $\tau_0,x_b,\kappa,\lambda$ for which \Cref{alg:mvc} achieves polylogarithmic regret under CI.

Specifically, the learner is supplied with constants $\alpha_0,\alpha_1,\vartheta,\nu,\nu_0>0$ such that
\[
 \|Q\|,\|R\|\le\alpha_1,\quad \|Q^{-1}\|\le\alpha_0^{-1},\quad
 \|A_\star\|,\|B_\star\|\le\vartheta,\quad
 J_\star\le\nu,\qquad J(K_0)\le\nu_0. 
\]
Moreover, it is given conservative lower bounds on $\lambda_{\min}(S_\star)$ and on $\sigma_{\min}^{+}(D_\star)$, the smallest nonzero singular value of the linear map $D_\star h=\Delta A_h+\Delta B_hK_\star$, which gives the change of the closed loop under the fixed gain $K_\star$ when the parameter moves from $\theta_\star$ along $h\in\mathbb R^p$ (the latter bound is not needed if $D_\star=0$).
The resulting values of the algorithmic parameters are 
\begin{equation}\label{eq:algorithmic_parameters}
    \tau_0=\lceil C_\tau\log(5T)\rceil,\quad x_b=C_x\log(5T),\quad \lambda=C_\lambda\log(5T),
\end{equation}
where $\kappa$ and the constants $C_\tau,C_x,C_\lambda$ depend on these bounds but not on $T$ (see \Cref{remark:constants}).

As in~\citet{casselLogarithmicRegretLearning2020} and~\citet{lee2024nonasymptoticregretanalysisadaptive}, the warm-up provides an initial estimate that is accurate enough for the CE gains to be well defined and stabilizing, which is needed to ensure that, with high probability, the algorithm never falls back to $K_0$. %
Unlike an explore-then-commit scheme that fixes its controller after exploring, \Cref{alg:mvc} keeps updating its estimate after the warm-up using the data generated by CE, so the warm-up only needs to reach a fixed accuracy that does not depend on $T$. This accuracy must, however, be reached with a confidence that increases with $T$.
An inaccurate warm-up estimate may cause the algorithm to fall back to $K_0$ for the remaining rounds, incurring regret of order $T$.
A fixed failure probability would therefore lead to a contribution of order $T$ to expected regret.
Our analysis chooses a failure probability $\delta=T^{-2}$, which keeps the expected contribution of failures bounded (Appendix~\ref{appendix:proof_mvc_regret_first_term}).
Since a fixed accuracy can be reached with failure probability $\delta$ using $O(\log(1/\delta))$ samples, this choice gives a warm-up length $\tau_0=O(\log T)$. %

\begin{algorithm}
\caption{Certainty Equivalence for a CI system}
\begin{algorithmic}[1]\label{alg:mvc}
\Require horizon $T$, parameters $\tau_0,x_b,\kappa,\lambda$,
         a stabilizing gain $K_0$, and the affine class $\mathcal C$
\State \textbf{initialize:} $N_T=\left\lceil (T-\tau_0+1)/\tau_0\right\rceil$, $\tau_i\gets\min\{\tau_0(i+1),T+1\}$ for $0\le i\le N_T$, and $\mathrm{aborted}\gets\mathrm{false}$
\For{$t = 1, \ldots, \min\{\tau_0-1,T\}$} 
    \State play $u_t = K_0 x_t+\eta_t$, where $\eta_t\overset{\text{i.i.d.}}{\sim}\mathcal N(0,\sigma^2I_m)$
\EndFor

\For{epoch $i = 0, \ldots, N_T-1$} 
    \State compute $\hat\Theta_{\tau_i}$ by \eqref{eq:constrained_ridge}
    \If{no $\hat P_{\tau_i}\succ0$ solves the estimated DARE with $\hat S_{\tau_i}\succ0$ and stable estimated closed loop, or $\|K(\hat\Theta_{\tau_i})\|>\kappa$}
        \State $\mathrm{aborted}\gets\mathrm{true}$
    \EndIf
    \For{$t = \tau_i, \ldots, \tau_{i+1}-1$}
        \State \textbf{if} $\|x_t\|^2>x_b$ \textbf{then} $\mathrm{aborted}\gets\mathrm{true}$
        \State play $u_t=K_0x_t$ if $\mathrm{aborted}$, else $u_t=K_{\tau_i}x_t$ with $K_{\tau_i}=K(\hat\Theta_{\tau_i})$
    \EndFor
\EndFor
\end{algorithmic}
\end{algorithm}

\subsubsection{Regret guarantee}
\begin{theorem}\label{thm:mvc_logregret}
Let \Cref{assumption:riccati} hold, and suppose CI holds at
$\theta_\star$.
Run \Cref{alg:mvc} with the parameters of \Cref{eq:algorithmic_parameters}.
Then, for every $T\geq1$,
\[
 \mathbb E[R_T]\leq C\log^2(5T),
\]
where $C<\infty$ depends on the problem parameters
and the supplied bounds, but not on $T$.
\end{theorem}

\paragraph{Proof sketch.} The full proof, including an explicit expression for $C$, is given in Appendix~\ref{appendix:proof_mvc_regret_first_term}.
We first bound the expected regret incurred during the warm-up by $O(\tau_0)=O(\log T)$.
We then analyze the remaining rounds separately on a success event and its complement.
On the success event, the warm-up produces an estimate sufficiently close to the true model that its CE gain is well defined and stabilizing.
State bounds and estimation confidence bounds ensure that subsequent estimates remain in this neighborhood, so the algorithm never aborts.
CI then bounds the suboptimality of each estimated controller by the closed-loop prediction error of the estimated model under that controller.
A self-normalized confidence bound~\citep{pmlr-v19-abbasi-yadkori11a} and the elliptical-potential lemma~\citep{AbbasiYadkori2011ImprovedAF} bound the cumulative prediction error, each contributing one logarithmic factor, which yields the $O(\log^2 T)$ term.
Finally, the choices of $\tau_0,x_b,\lambda$ ensure that the failure event has probability at most $T^{-2}$.
The abort mechanism controls the cost incurred on this event by switching to $K_0$, yielding an $O(1)$ contribution to expected regret.

\paragraph{Knowledge of the horizon.}
Similarly to~\cite{casselLogarithmicRegretLearning2020}, \Cref{alg:mvc} requires knowledge of the horizon $T$, which is used to ensure that the probability of the failure event is at most $T^{-2}$ and that the states collected in every epoch are sufficiently informative. We expect that this requirement can be removed, e.g., by 
allowing the algorithm to resume CE even after using $K_0$~\citep{jedra2021minimalexpectedregretlinear}, and by using epochs that only end once sufficiently informative data have been collected~\citep{bartosOptimisticOnlineLQR2026}. We leave this extension for future work.

\paragraph{Verifying CI.}
\Cref{thm:mvc_logregret} requires CI at the unknown parameter $\theta_\star$, which in general cannot be checked before learning. In some problem classes, however, CI follows from the prior knowledge available to the learner. For example, if $B_\star$ is known, as in~\citet{casselLogarithmicRegretLearning2020} and~\citet{jedra2021minimalexpectedregretlinear}, one can show that $\mathcal{I}_\star\succ 0$, so CI holds for every $R\succeq0$. Their setting with known $A_\star$ also satisfies CI under their additional condition $K_\star K_\star^\top\succ0$. %
In both settings, however, all parameters are identifiable under optimal feedback. We show in \Cref{sec:mvc} that, for MVC ($R=0$), CI holds whenever $\operatorname{range}(B_\star)$ is known, an interpretable condition on the model class under which up to $m^2$ parameter directions remain invisible under optimal feedback.

\paragraph{Comparison with local minimax guarantees.} 
\citet{liu2026controller} proves local minimax guarantees for $R\succ0$:  
assuming CI at a known nominal parameter, $O(\log T)$ regret bounds hold uniformly over neighborhoods
 of radius $O(T^{-1/4})$ around it, and over neighborhoods of sufficiently small fixed radius if the nominal state matrix is nonsingular,
 although in this case CI excludes nontrivial invisible directions.
Our algorithm instead requires no nominal model and achieves $O(\log^2 T)$ expected regret at every fixed unknown parameter satisfying CI, but its constant degrades as the true system approaches one where CI fails (\Cref{section:mvc_lowerbound}).
Hence, neither type of guarantee implies the other.

\subsection{Regret Lower Bound} \label{section:mvc_lowerbound}

We now show that the bound of \Cref{thm:mvc_logregret} cannot hold uniformly over the following family of systems, which includes a system where CI fails. Fix
$r\geq0$ and consider the family defined by the affine model
\begin{equation}\label{eq:lower_affine_family}
 A=\begin{bmatrix}0&1/4\\0&0\end{bmatrix},\qquad
 B_\theta=\begin{bmatrix}\theta\\1\end{bmatrix},\qquad
 |\theta|\leq0.05,
\end{equation}
with $Q=I_2$, $R=r$, $\sigma=1$, and $x_1=0$, where 
  $\theta$ is the unknown parameter and $K_0=0$ stabilizes every system in the family.
At $\theta=0$, the optimal controller applies zero input. Since $\theta$ appears only in the input matrix, observations under this controller contain no information about $\theta$. Nevertheless, controllers designed for nearby parameters incur a quadratic increase in cost on the $\theta=0$ system, so this invisible direction is not harmless and CI fails. For every $\theta\neq0$, the scalar Fisher information $\mathcal I_\theta(K(\theta))$ is positive, so there is no nonzero invisible parameter direction and CI holds. However, as $\theta$ approaches zero, the optimal controller uses progressively smaller inputs, so the data become less informative about $\theta$. Consequently, the constant in \Cref{thm:mvc_logregret} deteriorates even though CI continues to hold at every nonzero $\theta$ (\Cref{appendix:mvc_lower_bound_proof}).
 The following theorem shows that this cannot be avoided.
\begin{theorem}\label{thm:lower_bound}
For every fixed $r\geq0$ and every $T\geq10^4$, over the above family of systems
\[
 \inf_\pi\sup_{|\theta|\leq0.05}
 \mathbb E_\theta^\pi[R_T]
 \geq \frac{3.6\cdot10^{-4}}{1+r}\sqrt T-2.4,
\]
where the infimum is over all causal, possibly randomized policies.
\end{theorem}

\paragraph{Proof sketch.} \looseness-1 The full proof is given in Appendix~\ref{appendix:mvc_lower_bound_proof}.
We compare two systems with $\theta=\pm\Delta$, where $\Delta=T^{-1/4}/2$.
To keep regret small, the learner must choose inputs close to optimal feedback.
In both systems, these inputs are small, so the resulting observations provide too little information to distinguish the two parameter values reliably.
Distinguishing them matters because their optimal controllers prescribe opposite actions.
These actions differ by order $\Delta$ whenever the second state component is bounded away from zero, which Gaussian noise ensures with constant probability.
The learner therefore cannot stay close to both optimal controllers.
Averaging over the two systems, the resulting control errors yield a regret lower bound of order $T\Delta^2$, and hence $\sqrt T$, even when $R=0$.

\Cref{thm:lower_bound} gives a worst-case regret lower bound of order $\sqrt T$ over this family, which includes the system $\theta=0$, where CI fails.
The bound remains valid when $R=0$, so the absence of an input penalty alone does not guarantee polylogarithmic minimax regret.
Since both systems used in the proof have $\theta\neq0$, the same lower bound holds over the subfamily obtained by removing $\theta=0$, in which every system satisfies CI. This is consistent with \Cref{thm:mvc_logregret}, since its instance-specific constant $C$ is finite for every fixed $\theta\neq0$ but not uniformly bounded as $\theta$ approaches zero.

\section{Minimum Variance Control}\label{sec:mvc}
We consider MVC in the state-space setting of \Cref{sec:pb_setting}, with no input penalty ($R=0$). The objective is therefore to minimize the stationary weighted state variance $\mathbb E[x^\top Qx]$.
\begin{corollary}[Minimum variance control]\label{cor:mvc_known_range}
Let $R=0$, suppose \Cref{assumption:riccati} holds, and consider the model class
 $\mathcal C=\{[A\ B]: A\in\mathbb R^{n\times n},\ B\in\mathcal B\}$, where $\mathcal B\ni B_\star$ is a known affine set. 
If $\operatorname{range}(B)\subseteq\operatorname{range}(B_\star)$ for every $B\in\mathcal B$,
 then CI holds at $\theta_\star$, and the regret bound of \Cref{thm:mvc_logregret} applies.
  Conversely, if $\operatorname{rank}[A_\star\ B_\star]=n$, CI holds only if this range condition holds.
\end{corollary}

The proof is given in \Cref{appendix:mvc_proof}. 
The range condition means that the learner knows $\operatorname{range}(B_\star)$, the subspace through which inputs enter the state equation.
Let $U$ have orthonormal columns spanning this subspace.
The largest affine set satisfying the range condition is
 $\mathcal B=\{U\Gamma:\Gamma\in\mathbb R^{m\times m}\}$.
In this class, $\ker\mathcal I_\star$ has dimension $m^2$.
\Cref{cor:mvc_known_range} guarantees CI, so all these invisible directions are harmless.
When $A_\star$ is nonsingular, this attains the bound of \Cref{lemma:degenerate_cost_dimension}, whereas CI would fail for this model class under any positive definite input penalty $R\succ0$.

Knowing $\operatorname{range}(B_\star)$ requires less prior information than knowing the full matrix $B_\star$, as assumed in the corresponding result of \citet{casselLogarithmicRegretLearning2020} for $R\succ0$.
For square, nonsingular input matrices, as assumed in \citeauthor{LAI86}'s (\citeyear{LAI86}) stochastic regression model, the range is $\mathbb R^n$, so the range condition holds automatically, whereas our corollary also covers underactuated systems with $m<n$, provided $\operatorname{range}(B_\star)$ is known.
Finally, note that in the lower-bound family of \Cref{section:mvc_lowerbound}, $A_\star$ is known but $\operatorname{range}(B_\star)$ is not. There, the constant of \Cref{thm:mvc_logregret} degrades because $\sigma_{\min}^{+}(D_\star)\to0$ as $\theta\to0$. In contrast, $\sigma_{\min}^{+}(D_\star)\geq1$ in the class of \Cref{cor:mvc_known_range} (\Cref{appendix:mvc_proof}) and $\lambda_{\min}(S_\star)\geq\alpha_0\sigma_{\min}(B_\star)^2$. Hence, the constant $C$ of \Cref{thm:mvc_logregret} is uniformly bounded over all systems in this class that satisfy the norm and cost bounds of \Cref{sec:prior_information} and have $\sigma_{\min}(B_\star)$ bounded away from zero (\Cref{remark:mvc_constants}).

\subsection{Cart--Pole Case Study}\label{section:cartpole_case_study}
\begin{figure}
    \centering
    \includegraphics{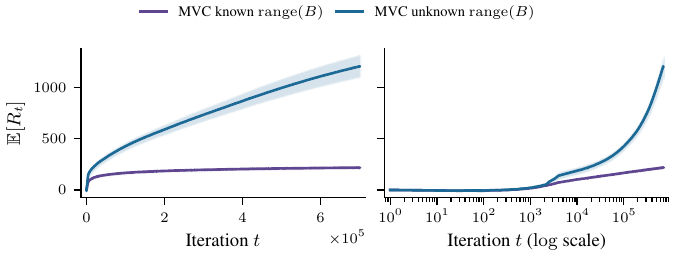}
    \caption{Expected regret of \Cref{alg:mvc}, estimated via \eqref{eq:performance_difference_appendix}, for minimum variance control of the linearized cart--pole, with and without knowledge of $\operatorname{range}(B_d)$. We report the mean and the standard error over $100$ seeds. We use $h=0.05$, $M_c=1$, $M_p=0.1$, $\ell=0.5$, $g=9.81$, $g_u=20$, $Q=I_4$, $\sigma=0.05$, and $T=7\cdot10^5$; since the constants of \Cref{remark:constants} are conservative, we set $\tau_0=2000$ and choose the norm thresholds heuristically (they were never triggered in the reported runs).}
    \label{fig:cartpole_linear}
\end{figure}

We consider a case study of learning to balance a pole upright on a cart.
The cart has mass $M_c>0$, and the pole is modeled as a point mass $M_p>0$ attached to the cart by a massless rigid rod of length $\ell>0$.
The controller applies a horizontal force $F=g_u u$ to the cart, where $u\in\mathbb R$ is the controller command and $g_u>0$ is the actuator gain.
We use the state $x=(p,\phi,\dot p,\dot\phi)^\top$, where $p$ is the cart position and $\phi$ is the pole's angular deviation from upright.
We linearize the equations of motion in \citet[Ch.~3]{tedrake2024underactuated} around the unstable upright equilibrium $(x,u)=(0,0)$, with $g$ denoting gravitational acceleration.
Discretizing the resulting dynamics using the forward-Euler method with sampling interval $h>0$ gives
\begin{equation*}
 \begin{gathered}A_d=I_4+h\begin{bmatrix}
 0&0&1&0\\
 0&0&0&1\\
 0&\frac{M_pg}{M_c}&0&0\\
 0&\frac{(M_c+M_p)g}{M_c\ell}&0&0
 \end{bmatrix},\quad
 B_d=h\frac{g_u}{M_c}
 \begin{bmatrix}0\\0\\1\\1/\ell\end{bmatrix}.
 \end{gathered}
\end{equation*}
\looseness=-1 Notably, the pole length, which can be easily measured in practice, suffices to determine
 $\operatorname{range}(B_d)$, whereas identifying $B_d$ would additionally 
require the scale $hg_u/M_c$.
Since $A_d$ is nonsingular for small $h$, \Cref{cor:mvc_known_range} implies that CI holds for MVC if and only if
 the learner knows $\operatorname{range}(B_d)$.
\Cref{fig:cartpole_linear} compares two learners: one knows $\operatorname{range}(B_d)$ and estimates only the scaling, while the other estimates the last two entries of $B_d$ independently.
With the range known, the observed regret is consistent with the polylogarithmic upper bound.
Without this information, regret grows visibly faster.

\section{Conclusion}
This work studied when certainty equivalence achieves polylogarithmic regret in online LQR with possibly singular input penalties.
We extended controller identifiability to this setting and showed that it suffices for $O(\log^2 T)$ expected regret.
Our results clarify how singular input penalties can make parameter directions that are unidentifiable under optimal feedback harmless for control.
In MVC, knowing $\operatorname{range}(B_\star)$ suffices for controller identifiability, allowing polylogarithmic regret without identifying all system parameters.
A complementary lower bound shows that the absence of an input penalty alone does not guarantee polylogarithmic minimax regret.

An interesting direction for future work is to carry these insights over to nonlinear control by extending the analysis of \citet{wang2025logarithmicregretnonlinearcontrol}.
In particular, one could investigate whether their persistent-excitation condition under the optimal policy can be relaxed to a nonlinear analogue of CI.

\subsection*{Generative AI Use}
Generative AI tools have only been used in this work for text polishing, code assistance, and checking mathematical proofs. All ideas, derivations, and conclusions presented in this work remain the sole responsibility of the authors.

\acks{Bruce D. Lee was supported by an ETH AI Center Postdoctoral Fellowship. This work was also supported in part by NCCR Automation, a National Centre of Competence in Research, funded by the Swiss National Science Foundation (grant number 51NF40\_225155).}

\bibliography{references}

\clearpage

\appendix
\crefalias{section}{appendix}\crefalias{subsection}{subappendix}
\section{Technical Lemmas}\label{appendix:a}
\label{sec:a1_general}
This appendix collects properties of the Riccati equation for possibly singular $R$. We retain the notation from \Cref{sec:pb_setting}  and impose
\Cref{assumption:riccati}, with $Q\succeq\alpha_0I\succ0$, $R\succeq0$, and
$\max\{\|Q\|,\|R\|\}\leq\alpha_1$ as in \Cref{sec:prior_information}. Lemma, equation, and section numbers of~\cite{mania2019extended} refer to its extended version, arXiv:1902.07826.

\subsection{Riccati Well-Posedness}
The following proposition collects well-known properties of the discrete algebraic Riccati equation, under the conditions of \Cref{assumption:riccati}. A short proof is provided below for completeness.
\begin{proposition}\label{thm:riccati_solution}\label{lemma:general_riccati_identification}
  Suppose that $(A,B)$ is stabilizable, $Q\succ0$, $R\succeq0$ and
  $\ker R\cap\ker B=\{0\}$. Then the DARE
  \begin{equation}\label{eq:singular_dare}
      P=Q+A^\top PA-A^\top PB(R+B^\top PB)^{-1}B^\top PA
  \end{equation}
  has a unique positive definite solution $P$, and, for every $x_0\in\mathbb R^n$,
  \begin{equation}\label{eq:value_function}
   x_0^\top Px_0=\inf_{u_0,u_1,\dots}\sum_{t\ge0}\bigl(x_t^\top Qx_t+u_t^\top Ru_t\bigr),
  \end{equation}
  where $x_{t+1}=Ax_t+Bu_t$. With $K=-(R+B^\top PB)^{-1}B^\top PA$, the matrix $A+BK$ is Schur stable, $u_t=Kx_t$ is the unique stationary linear policy attaining the infimum in \eqref{eq:value_function} from every initial state, and $K$ is the unique minimizer of the average cost $J_{[A\ B]}$ in \eqref{eq:cost_definition}. Conversely, if $Q\succ0$, $R\succeq0$, $\ker R\cap\ker B=\{0\}$, and \eqref{eq:singular_dare} has a positive definite solution, then $(A,B)$ is stabilizable.
\end{proposition}
\begin{proof}
Stabilizability makes the cost in \eqref{eq:value_function} finite for every $x_0$, so \citet[Theorem~2.1]{ferrante2012generaliseddiscretealgebraicriccati} gives a solution $P\succeq0$ of the generalized DARE satisfying \eqref{eq:value_function}, and $P\succeq Q\succ0$ since the cost is at least $x_0^\top Qx_0$. By the kernel assumption, $R+B^\top PB\succ0$, so $P$ solves \eqref{eq:singular_dare} and, by the same theorem, $u_t=Kx_t$ is the unique optimal input.
Since $Q\succ0$, Lyapunov's theorem shows that every positive definite solution of \eqref{eq:singular_dare} is stabilizing, which gives the Schur stability of $A+BK$ and the converse; uniqueness then follows from~\citet[Proposition~1]{articleIonescu}.

For the average cost, quadratic completion gives, with $S=R+B^\top PB$,
\begin{equation}\label{eq:riccati_completion}
 x^\top Qx+u^\top Ru+(Ax+Bu)^\top P(Ax+Bu)-x^\top Px
 =(u-Kx)^\top S(u-Kx).
\end{equation}
For a gain $G$ stabilizing $(A,B)$, the Lyapunov and covariance equations give, with $J=J_{[A\ B]}$,
\begin{equation}\label{eq:a1_policy_evaluation}
\begin{split}
 P_G&=Q+G^\top RG+(A+BG)^\top P_G(A+BG),\\
 J(G)&=\sigma^2\operatorname{Tr}(P_G),\qquad P_G\succ0,
\end{split}
\end{equation}
and taking stationary expectations in \eqref{eq:riccati_completion} with $u=Gx$ yields
\begin{equation}\label{eq:a1_performance_difference_stationary}
 J(G)-\sigma^2\operatorname{Tr}(P)
 =\operatorname{Tr}\!\bigl(\Sigma_G(G-K)^\top
              (R+B^\top PB)(G-K)\bigr),
\end{equation}
where $\Sigma_G=\sigma^2\sum_{j\geq0}(A+BG)^j((A+BG)^j)^\top\succeq\sigma^2I$ is the stationary state covariance.
This is the identity of \citet[Lemma~3, Eq.~(8)]{mania2019extended}, which assumes $R\succ0$; the derivation above uses only $R+B^\top PB\succ0$. Hence $J(G)>J(K)$ for every stabilizing $G\neq K$. If instead $\rho(A+BG)\geq1$, then $\|(A+BG)^j\|_F\geq1$ for all $j$, so $\mathbb E\|x_t\|^2\geq\sigma^2(t-1)$ and, since $Q\succ0$, the average cost is infinite.
\end{proof}
For any causal policy with finite second moments, taking
expectations in \eqref{eq:riccati_completion} at $(A_\star,B_\star)$ and telescoping gives
\begin{equation}\label{eq:performance_difference_appendix}
\mathbb E[R_T]
=\mathbb E\!\left[
  \sum_{t=1}^T
  (u_t-K_\star x_t)^\top S_\star(u_t-K_\star x_t)
\right]
+x_1^\top P_\star x_1
  -\mathbb E[x_{T+1}^\top P_\star x_{T+1}].
\end{equation}
Here independence and centering of $w_t$ contribute
$\sigma^2\operatorname{Tr}(P_\star)=J_\star$ per step.
This is the matrix counterpart of
\citet[Lemma~31]
{casselLogarithmicRegretLearning2020}.

\subsection{Lyapunov Inverse and Matrix Powers}
\begin{lemma}[cf.~{\citealp[Lemmas~4--5]{mania2019extended}}]
\label{lemma:general_lyapunov_inverse}\label{lemma:bound_inverse_operator}
Let $L_\star$ be Schur stable, fix
$\gamma_\star\in(\rho(L_\star),1)$, and put
\[
 \tau_\star\geq\sup_{j\geq0}\|L_\star^j\|\gamma_\star^{-j},
 \qquad L_{\mathcal T}=\frac{\tau_\star^2}{1-\gamma_\star^2}.
\]
The operator $\mathcal T(X)=X-L_\star^\top XL_\star$ is invertible and
\begin{equation}
 \mathcal T^{-1}(Y)=\sum_{j=0}^{\infty}(L_\star^\top)^jYL_\star^j,
 \qquad
 \|\mathcal T^{-1}(Y)\|_F\leq L_{\mathcal T}\|Y\|_F.
 \label{eq:general_lyapunov_inverse}
\end{equation}
Moreover, for every square matrix $E$ and integer $j\geq0$,
\begin{equation}
 \|(L_\star+E)^j\|
 \leq\tau_\star(\gamma_\star+\tau_\star\|E\|)^j.
 \label{eq:general_perturbed_powers}
\end{equation}
\end{lemma}
\begin{proof}
The series representation and the spectral-norm bound are in \citet[App.~C, proof of Lemma~4]{mania2019extended},
 and the power estimate is \citet[Lemma~5]{mania2019extended} with $M=L_\star$, $\rho=\gamma_\star$, and $\Delta=E$;
  finiteness of $\tau_\star$ follows from Gelfand's formula \citep[Eq.~(3)]{mania2019extended}. 
  The Frobenius bound follows from the same series, since $\|(L_\star^\top)^jYL_\star^j\|_F\leq\tau_\star^2\gamma_\star^{2j}\|Y\|_F$. 
  These statements involve only $L_\star$, so they do not depend on the standing assumption $R\succ0$ of~\cite{mania2019extended}.
\end{proof}
For the optimal closed loop, the Bellman equation $P_\star=Q+K_\star^\top RK_\star+L_\star^\top P_\star L_\star$ gives $L_\star^\top P_\star L_\star\preceq(1-\alpha_0/\|P_\star\|)P_\star$, hence $\|L_\star^j\|\leq\sqrt{\|P_\star\|/\lambda_{\min}(P_\star)}\,(1-\alpha_0/\|P_\star\|)^{j/2}$. Since $\sigma^2\|P_\star\|\leq J_\star\leq\nu$ and $P_\star\succeq Q\succeq\alpha_0I$, the computable pair
\[
 \tau_\star=\sqrt{\frac{\nu}{\alpha_0\sigma^2}},\qquad
 \gamma_\star=\sqrt{1-\frac{\alpha_0\sigma^2}{2\nu}}
\]
is admissible, and $\|L_\star\|\leq\tau_\star$.

\subsection{Riccati Sensitivity}
\label{sec:general_riccati_sensitivity}
In this section, we introduce some notation and derive several technical preliminary identities that will be used in the proof of the Riccati sensitivity bound in \Cref{lemma:general_riccati_sensitivity}.
For a symmetric matrix $P$, define
\begin{align}
 S_{A,B}(P)&=R+B^\top PB,\notag\\
 K_{A,B}(P)&=-S_{A,B}(P)^{-1}B^\top PA,\notag\\
 L_{A,B}(P)&=A+BK_{A,B}(P),\notag\\
 \mathcal F_{A,B}(P)
 &=Q+A^\top PA-A^\top PB S_{A,B}(P)^{-1}B^\top PA,
 \label{eq:general_riccati_map}
\end{align}
whenever $S_{A,B}(P)\succ0$.  Under the standing assumptions,
\Cref{thm:riccati_solution} supplies the reference solution
$P_\star=\mathcal F_{A_\star,B_\star}(P_\star)\succ0$ and a stabilizing gain,
with
\[
 S_\star=R+B_\star^\top P_\star B_\star\succ0,\qquad
 K_\star=K_{A_\star,B_\star}(P_\star),\qquad L_\star=A_\star+B_\star K_\star.
\]

For $P\succ0$ and $S=R+B^\top PB\succ0$, write
$K=K_{A,B}(P)$, $L=A+BK$, $W=BS^{-1}B^\top$, and
$\ell=\lambda_{\min}(P)$.
The matrix $C=P^{1/2}BS^{-1/2}$ satisfies $C^\top C\preceq I$
because $S\succeq B^\top PB$. Using $CC^\top\preceq I$ and
$I-WP=P^{-1/2}(I-CC^\top)P^{1/2}$ therefore gives
\begin{equation}\label{eq:general_elimination_bounds}
\begin{gathered}
 0\preceq P^{1/2}WP^{1/2}\preceq I,\\
 \|W\|\leq\ell^{-1},\qquad \|I-WP\|\leq\sqrt{\|P\|/\ell},\\
 \|BS^{-1/2}\|\leq\ell^{-1/2},\qquad
 \|S^{-1/2}B^\top P\|\leq\sqrt{\|P\|}.
\end{gathered}
\end{equation}
For symmetric $X$ with $\|X\|<\ell$, attach the subscript $X$
to quantities evaluated at $P+X$.
Since $P+X\succeq(1-\|X\|/\ell)P$, subtracting the stationarity
equations  $S_XK_X=-B^\top(P+X)A$ and
$S_XK=SK+B^\top XBK$, and applying the matrix inversion lemma 
\citep[Appendix~C.4.3]{boyd2004convex} to
$I+(BS^{-1})(B^\top X)$ gives
\begin{equation}\label{eq:general_fixed_model_identities}
\begin{gathered}
 S_X\succeq(1-\|X\|/\ell)S\succ0,\qquad
 K_X-K=-S_X^{-1}B^\top XL,\\
 L_X=(I-BS_X^{-1}B^\top X)L=(I+WX)^{-1}L.
\end{gathered}
\end{equation}
Here $\|WX\|<1$ ensures invertibility.
Quadratic completion, valid for every gain $G$, gives
\begin{equation}\label{eq:general_completion_arbitrary_P}
 Q+G^\top RG+(A+BG)^\top P(A+BG)
 =\mathcal F_{A,B}(P)+(G-K)^\top S(G-K).
\end{equation}
For symmetric $Z$, the Riccati derivatives are
\begin{align}
 D_P\mathcal F_{A,B}(P)[Z]&=L^\top ZL,\notag\\
 D_{A,B}\mathcal F_{A,B}(P)[\Delta A,\Delta B]
 &=(\Delta A+\Delta BK)^\top PL
      +L^\top P(\Delta A+\Delta BK).
 \label{eq:general_riccati_derivatives}
\end{align}
The first is \citet[Eq.~(1.2)]{guo1998newton}, whose domain only
requires $S$ invertible. For the second, apply
\eqref{eq:general_completion_arbitrary_P} at
$(A+t\Delta A,B+t\Delta B)$ with $G=K$ fixed.
The gain $K_t=K_{A_t,B_t}(P)$ and input curvature $S_t=R+B_t^\top PB_t$ are smooth near zero because $S\succ 0$,
so $(K-K_t)^\top S_t(K-K_t)=O(t^2)$ because $K_t-K=O(t)$ and $S_t$ is bounded.
Differentiating the remaining quadratic gives the stated formula.

Write $\Delta A=\hat A-A_\star$, $\Delta B=\hat B-B_\star$,
$E=\Delta A+\Delta BK_\star$, $M=\Delta B^\top P_\star L_\star$,
$S_0=R+\hat B^\top P_\star\hat B$, and
$G_0=M+\hat B^\top P_\star E$.
Put $\ell_\star=\lambda_{\min}(P_\star)$ and
$s_\star=\lambda_{\min}(S_\star)$.
The reverse triangle inequality applied to the norm of 
$[R^{1/2}v;P_\star^{1/2}\hat Bv]$ gives
\begin{equation}\label{eq:general_reference_curvature}
 S_0\succeq(\sqrt{s_\star}-\sqrt{\|P_\star\|}\|\Delta B\|)^2I\succ0\qquad
 \text{if }\sqrt{\|P_\star\|}\|\Delta B\|<\sqrt{s_\star}.
\end{equation}
Reference stationarity yields
\[
RK_\star+\hat B^\top P_\star(L_\star+E)
=\Delta B^\top P_\star L_\star+\hat B^\top P_\star E
=G_0,
\]
hence
\begin{equation}\label{eq:general_reference_gain}
 K_{\hat A,\hat B}(P_\star)-K_\star=-S_0^{-1}G_0.
\end{equation}
Substituting in \eqref{eq:general_completion_arbitrary_P} and
subtracting the reference Bellman identity $P_\star=Q+K_\star^\top RK_\star+L_\star^\top P_\star L_\star$ gives
\begin{align}
 \mathcal F_{\hat A,\hat B}(P_\star)-P_\star
 &=L_\star^\top P_\star E+E^\top P_\star L_\star
      +E^\top P_\star E-G_0^\top S_0^{-1}G_0,
 \label{eq:general_reference_residual_identity}\\
 L_0-L_\star
 &=(I-\hat B S_0^{-1}\hat B^\top P_\star)E-\hat B S_0^{-1}M,
 \label{eq:general_reference_loop_identity}
\end{align}
where $L_0=L_{\hat A,\hat B}(P_\star)$ and we used $L_0=L_\star+E+\hat B(K_{\hat A,\hat B}(P_\star)-K_\star)$.
For any $s_0>0$ with $S_0\succeq s_0I$,
\eqref{eq:general_elimination_bounds} implies
\begin{equation}\label{eq:general_scaled_residual_bound}
 \|S_0^{-1/2}G_0\|_F
 \leq\sqrt{\|P_\star\|}\|E\|_F+\|M\|_F/\sqrt{s_0},\qquad
 \|M\|_F\leq\|P_\star L_\star\|\|\Delta B\|_F.
\end{equation}
Consequently, the residual identities give
\begin{align}
 \|\mathcal F_{\hat A,\hat B}(P_\star)-P_\star\|_F
 &\leq2\|P_\star L_\star\|\|E\|_F
   +\|P_\star\|\|E\|_F^2+\|S_0^{-1/2}G_0\|_F^2,
 \label{eq:general_reference_residual_bound}\\
 \|L_0-L_\star\|_F
 &\leq\sqrt{\|P_\star\|/\ell_\star}\|E\|_F
       +\frac{\|M\|_F}{\sqrt{\ell_\star s_0}}.
 \label{eq:general_reference_loop_bound}
\end{align}

The following perturbation argument follows the Lyapunov-inverse fixed-point
method in \citet[Proposition~2]{mania2019extended}.
Here the residual and derivative estimates use $R+B^\top PB\succ0$,
allowing $R$ to be singular.

\begin{lemma}
\label{lemma:general_riccati_sensitivity}
Let $\gamma_\star,\tau_\star,L_{\mathcal T}$ be as in \Cref{lemma:general_lyapunov_inverse}, and let $L_P$ and the radius $r>0$ be given by \eqref{eq:general_closed_form_constants}--\eqref{eq:general_closed_form_radius} below; they depend only on $\|P_\star\|$, $\ell_\star$, $s_\star$, $\|K_\star\|$, $\|L_\star\|$, $\|P_\star L_\star\|$, and $L_{\mathcal T}$. If $\varepsilon=\|[\hat A-A_\star\ \hat B-B_\star]\|_F\leq r$, then $(\hat A,\hat B)$ has a unique stabilizing Riccati solution $\hat P\succ0$, with $\hat S=R+\hat B^\top\hat P\hat B\succ0$, and
\[
 \|\hat P-P_\star\|_F\leq L_P\varepsilon;
\]
moreover, \eqref{eq:general_P_S_loop_bounds} below holds, with $\eta_r,s_r,\xi_r$ defined in \eqref{eq:general_sensitivity_constants}. These conclusions remain valid when $r$ is replaced by any smaller radius in all the quantities below. 
\end{lemma}

\begin{proof}
Put $c_E=\sqrt{1+\|K_\star\|^2}$, and define the following bounds, used to choose the radius:
\begin{align}
 L_P&=2L_{\mathcal T}(1+2\|P_\star L_\star\|c_E),\notag\\
 \bar N&=\|P_\star\|c_E^2+
       \left(\sqrt{\|P_\star\|}c_E+
              \frac{2\|P_\star L_\star\|}{\sqrt{s_\star}}\right)^2,\notag\\
 \bar\xi&=2\left(\sqrt{\frac{\|P_\star\|}{\ell_\star}}c_E+
       \frac{2\|P_\star L_\star\|}{\sqrt{\ell_\star s_\star}}+
       \frac{\|L_\star\|L_P}{\ell_\star}\right).
 \label{eq:general_closed_form_constants}
\end{align}
Set
\begin{equation}
 r=\min\left\{
 \frac{\ell_\star}{2L_P},\quad
 \frac{\sqrt{s_\star}}{2\sqrt{\|P_\star\|}},\quad
 \frac{1}{\bar N},\quad
 \frac{1}{2L_{\mathcal T}\bar\xi
   (\sqrt{\|L_\star\|^2+(2L_{\mathcal T})^{-1}}+\|L_\star\|)}
 \right\}.
 \label{eq:general_closed_form_radius}
\end{equation}
We also define
\begin{align}
 \eta_r&=L_Pr/\ell_\star,\qquad
 s_r=(\sqrt{s_\star}-\sqrt{\|P_\star\|}r)^2,\notag\\
 \xi_r&=\frac{1}{1-\eta_r}\left(
       \sqrt{\frac{\|P_\star\|}{\ell_\star}}c_E+
       \frac{\|P_\star L_\star\|}{\sqrt{\ell_\star s_r}}+
       \frac{\|L_\star\|L_P}{\ell_\star}\right),\notag\\
 \chi_r&=L_{\mathcal T}(2\|L_\star\|\xi_rr+\xi_r^2r^2).
 \label{eq:general_sensitivity_constants}
\end{align}
For a perturbed model with $\varepsilon\leq r$, one has $\|\Delta B\|_F\leq\varepsilon$ and $\|E\|_F\leq c_E\varepsilon$.
The radius restrictions immediately give $\eta_r\leq1/2$, $s_r\geq s_\star/4$,
and $\xi_r\leq\bar\xi$.
The last restriction is equivalent to
$\bar\xi r\leq\sqrt{\|L_\star\|^2+(2L_{\mathcal T})^{-1}}-\|L_\star\|$,
which implies $\chi_r\leq1/2$.
Since $\bar N r\leq1$ and
$L_P/2=L_{\mathcal T}(1+2\|P_\star L_\star\|c_E)$, we obtain
the certificate
\begin{equation}\label{eq:general_sensitivity_certificate}
 \eta_r,\chi_r\leq\tfrac12,\qquad
 \sqrt{\|P_\star\|}r\leq\tfrac12\sqrt{s_\star},\qquad
 L_{\mathcal T}(2\|P_\star L_\star\|c_E+\bar N r)
 \leq(1-\chi_r)L_P.
\end{equation}
All denominators in \eqref{eq:general_closed_form_radius}
are positive, so $r>0$. The case $\varepsilon=0$ is immediate;
henceforth suppose $0<\varepsilon\leq r$.

Consider the closed ball of symmetric matrices
$\mathcal B_\varepsilon=\{X:\|X\|_F\leq L_P\varepsilon\}$,
and write $P_X=P_\star+X$,
$L_X=L_{\hat A,\hat B}(P_X)$, and
$S_X=R+\hat B^\top P_X\hat B$.
By \eqref{eq:general_reference_curvature} and
\eqref{eq:general_fixed_model_identities}, throughout this ball
\[
 S_0\succeq s_rI,\qquad P_X\succeq(1-\eta_r)P_\star\succ0,\qquad
 S_X\succeq(1-\eta_r)S_0\succ0.
\]
Thus the Riccati map is defined on the entire ball.
Put $W_0=\hat B S_0^{-1}\hat B^\top$.
Since $\|W_0X\|\leq\eta_r<1$,
\eqref{eq:general_fixed_model_identities} gives
\[
 L_X-L_\star=(I+W_0X)^{-1}
                  (L_0-L_\star-W_0XL_\star).
\]
The Neumann-series bound and \eqref{eq:general_reference_loop_bound}
therefore imply
\begin{equation}\label{eq:general_loop_on_ball}
 \|L_X-L_\star\|_F
 \leq\frac{\|L_0-L_\star\|_F+
              \ell_\star^{-1}\|X\|_F\|L_\star\|}
             {1-\eta_r}
 \leq\xi_r\varepsilon.
\end{equation}

With $\mathcal T(X)=X-L_\star^\top XL_\star$, define
\begin{align}
 \mathcal R(X)&=\mathcal F_{\hat A,\hat B}(P_\star+X)-P_\star
                   -L_\star^\top XL_\star,\notag\\
 \Phi(X)&=\mathcal T^{-1}(\mathcal R(X)).
 \label{eq:general_quantitative_fixed_point}
\end{align}
A fixed point of $\Phi$ is precisely a DARE solution $P_\star+X$.
By \eqref{eq:general_riccati_derivatives},
$D\mathcal R(X)[Z]=L_X^\top ZL_X-L_\star^\top ZL_\star$.
Expanding around $L_\star$ and using
\eqref{eq:general_loop_on_ball} bounds its norm by
$(2\|L_\star\|\xi_r\varepsilon+\xi_r^2\varepsilon^2)\|Z\|_F$.
Integrating along segments in $\mathcal B_\varepsilon$ and using
\eqref{eq:general_lyapunov_inverse} gives
\begin{equation}\label{eq:general_contraction_bound}
\begin{split}
 \|\Phi(X)-\Phi(Y)\|_F
 &\leq L_{\mathcal T}(2\|L_\star\|\xi_r\varepsilon+
                        \xi_r^2\varepsilon^2)\|X-Y\|_F\\
 &\leq\chi_r\|X-Y\|_F.
\end{split}
\end{equation}
Since $s_r\geq s_\star/4$,
\eqref{eq:general_scaled_residual_bound}--\eqref{eq:general_reference_residual_bound}
bound the quadratic residual coefficient by $\bar N$.
Together with the certificate, this gives
\[
\begin{aligned}
 \|\Phi(X)\|_F&\leq\|\Phi(0)\|_F+\chi_r\|X\|_F\\
 &\leq L_{\mathcal T}(2\|P_\star L_\star\|c_E+\bar N\varepsilon)
                    \varepsilon+\chi_rL_P\varepsilon\\
 &\leq L_P\varepsilon.
\end{aligned}
\]
Hence $\Phi$ maps $\mathcal B_\varepsilon$ into itself and is a
contraction. Banach's theorem supplies a fixed point $X$.
The matrix $\hat P=P_\star+X\succ0$ solves the DARE with
$\hat S\succ0$, which also implies $\ker R\cap\ker\hat B=\{0\}$.
The converse in \Cref{lemma:general_riccati_identification}
therefore proves stability and uniqueness without assuming
stabilizability of the perturbed model.
The ball and closed-loop bounds give $\|\hat P-P_\star\|_F\leq L_P\varepsilon$ and
\begin{equation}\label{eq:general_P_S_loop_bounds}
 \hat P\succeq(1-\eta_r)P_\star,\qquad \hat S\succeq(1-\eta_r)s_rI,\qquad
 \|\hat L-L_\star\|_F\leq\xi_r\varepsilon.
\end{equation}
Finally, every restriction in \eqref{eq:general_closed_form_radius} is an upper bound on the radius, so the argument applies verbatim with $r$ replaced by any smaller radius.
\end{proof}

The following corollary shows that the Riccati solution and the gain are smooth on an explicit neighborhood of the true model, and that the gain is uniformly Lipschitz in that neighborhood.

\begin{corollary}
\label{cor:general_uniform_sensitivity}
Let $r$ be the radius of \Cref{lemma:general_riccati_sensitivity}
and consider the convex set of models
\[
 \mathcal N_r=\{(A,B):\|[A-A_\star\ B-B_\star]\|_F\leq r\}.
\]
The stabilizing Riccati solution and the CE gain are real-analytic functions of $(A,B)$, in particular twice continuously differentiable,
in a neighborhood of every model in $\mathcal N_r$.
For any two models $(A_1,B_1),(A_2,B_2)\in\mathcal N_r$ with CE gains $K_1,K_2$,
\begin{equation}\label{eq:general_two_model_bounds}
 \|K_1-K_2\|_F
 \leq L_{K,1}\|[A_1-A_2\ B_1-B_2]\|_F,
\end{equation}
where $L_{K,1}$ is given in \eqref{eq:general_uniform_constants} below.
\end{corollary}
\begin{proof}
Put $\bar k_r=\sqrt{(\|P_\star\|+L_Pr)/((1-\eta_r)s_r)}\,(\|A_\star\|+r)$ and define
\begin{align}
 L_{P,1}&=\frac{2L_{\mathcal T}}{1-\chi_r}
  (\|P_\star\|+L_Pr)(\|L_\star\|+\xi_rr)
 \sqrt{1+\bar k_r^2},\notag\\
 L_{K,1}&=\frac{(\|P_\star\|+L_Pr)(\|L_\star\|+\xi_rr)}{(1-\eta_r)s_r}
 +\sqrt{\frac{\|P_\star\|+L_Pr}{(1-\eta_r)s_r}}
    \sqrt{1+\bar k_r^2}
 +\frac{(\|L_\star\|+\xi_rr)L_{P,1}}
                  {(1-\eta_r)\sqrt{\ell_\star s_r}}.
 \label{eq:general_uniform_constants}
\end{align}

Fix $(A,B)\in\mathcal N_r$ and write $P,K,L,S$ for its Riccati
quantities.
Put $\mathcal T_L(Z)=Z-L^\top ZL$.  Then
$(\mathcal T-\mathcal T_L)(Z)=L^\top ZL-L_\star^\top ZL_\star
 $, and, as in the proof of
\eqref{eq:general_contraction_bound} with $\varepsilon=r$,
$\|\mathcal T^{-1}(\mathcal T-\mathcal T_L)\|_{F\to F}
\leq L_{\mathcal T}(2\|L_\star\|\xi_rr+\xi_r^2r^2)=\chi_r<1$.
Factoring
$\mathcal T_L=\mathcal T[I-\mathcal T^{-1}(\mathcal T-\mathcal T_L)]$
and applying the Neumann series yields
\begin{equation}\label{eq:general_uniform_lyapunov_inverse}
 \|\mathcal T_L^{-1}\|_{F\to F}
 \leq\frac{L_{\mathcal T}}{1-\chi_r}.
\end{equation}
Since $D_P(P-\mathcal F_{A,B}(P))=\mathcal T_L$, the
analytic implicit-function theorem gives a real-analytic local solution, since the Riccati map is rational wherever $S\succ0$. Positivity of $P,S$ persists locally, and
\Cref{thm:riccati_solution} identifies this solution uniquely.

Differentiating the DARE and $RK+B^\top P(A+BK)=0$ in any
direction $\Delta=(\Delta A,\Delta B)$ gives, using the chain rule and \eqref{eq:general_riccati_derivatives}
\begin{align}
 \mathcal T_L(DP[\Delta])
 &=(\Delta A+\Delta BK)^\top PL
 +L^\top P(\Delta A+\Delta BK),
 \label{eq:general_differentiated_DARE}\\
 S\,DK[\Delta]
 &=-\bigl(\Delta B^\top PL+B^\top DP[\Delta]L
                    +B^\top P(\Delta A+\Delta BK)\bigr).
 \label{eq:general_differentiated_stationarity}
\end{align}
These are \citet[Lemma~2]{liu2026controller} up to the sign convention $u=-Kx$ and require only $S\succ 0$ so they remain valid here; see also
\citet[Section~5.2]{kao2020differentiation}.

By \eqref{eq:general_elimination_bounds}, \Cref{lemma:general_riccati_sensitivity}, and $\|A\|\leq\|A_\star\|+r$, we have $\|K\|\leq\|S^{-1/2}\|\|S^{-1/2}B^\top P\|\|A\|\leq\bar k_r$. Since $\|\Delta A+\Delta BK\|_F
\leq\sqrt{1+\|K\|^2}\|[\Delta A\ \Delta B]\|_F$,
\eqref{eq:general_uniform_lyapunov_inverse} gives
\[
\begin{aligned}
 \|DP[\Delta]\|_F
 &\leq\frac{2L_{\mathcal T}\|P\|\|L\|}{1-\chi_r}
        \|\Delta A+\Delta BK\|_F\\
 &\leq L_{P,1}\|[\Delta A\ \Delta B]\|_F.
\end{aligned}
\]
The three terms of \eqref{eq:general_differentiated_stationarity}
and \eqref{eq:general_elimination_bounds} give
\[
 \|DK[\Delta]\|_F
 \leq\frac{\|P\|\|L\|}{\lambda_{\min}(S)}\|\Delta B\|_F
 +\frac{\|L\|\|DP[\Delta]\|_F}
        {\sqrt{\lambda_{\min}(P)\lambda_{\min}(S)}}
 +\sqrt{\frac{\|P\|}{\lambda_{\min}(S)}}
                         \|\Delta A+\Delta BK\|_F.
\]
Substituting the bounds from \Cref{lemma:general_riccati_sensitivity} bounds the gain derivative by
$L_{K,1}\|[\Delta A\ \Delta B]\|_F$. Since $\mathcal N_r$ is convex,
integrating the gain derivative along a segment proves
\eqref{eq:general_two_model_bounds}.
\end{proof}

\subsection{Certainty Equivalence and Strong Stability}
The argument follows \citet[Section~2.3]{mania2019extended},
using the semidefinite-cost sensitivity bound above.
Following~\cite{casselLogarithmicRegretLearning2020}, we say that $K$ is $(\kappa,\gamma)$-strongly stable if
$\|K\|\leq\kappa$ and $A_\star+B_\star K=HLH^{-1}$ for some
$H,L$ with $\|H\|\|H^{-1}\|\leq\kappa$ and $\|L\|\leq1-\gamma$.

\begin{lemma}
\label{lemma:a1_ce}\label{lemma:a1_strong_ce}
Let $r$ be the radius of \Cref{lemma:general_riccati_sensitivity}, let $C_K=L_{K,1}$ be the constant of \Cref{cor:general_uniform_sensitivity},
and put
\[
 r_{\mathrm{CE}}=\min\left\{r,
       \frac{1-\gamma_\star}{2\tau_\star\|B_\star\|C_K}\right\},\qquad
 C_J=\frac{\sigma^2\tau_\star^2\|S_\star\|C_K^2}
              {1-((1+\gamma_\star)/2)^2},
\]
omitting the second radius restriction if $\|B_\star\|C_K=0$.
If a model $(A,B)$ satisfies $\varepsilon=\|[A-A_\star\ B-B_\star]\|_F\leq r_{\mathrm{CE}}$, then its CE gain $K$ exists, stabilizes the true plant, and satisfies
\[
 \|K-K_\star\|_F\leq C_K\varepsilon,\qquad
 J(K)-J_\star\leq C_J\varepsilon^2.
\]
Moreover, for $0<\varepsilon_0\leq r_{\mathrm{CE}}$, $J_\star\leq\nu$,
and $\bar k_\star\geq\|K_\star\|$, every such gain with
$\varepsilon\leq\varepsilon_0$ is $(\kappa,\gamma)$-strongly stable, where
\[
 \kappa=\max\left\{
 \sqrt{\frac{\nu+C_J\varepsilon_0^2}{\alpha_0\sigma^2}},
 \bar k_\star+C_K\varepsilon_0,1\right\},
 \qquad \gamma=\frac1{2\kappa^2}.
\]
\end{lemma}
\begin{proof}
\Cref{lemma:general_riccati_sensitivity} supplies existence, and \eqref{eq:general_two_model_bounds} with $(A_\star,B_\star)\in\mathcal N_r$ gives the gain bound.
The radius restriction and \eqref{eq:general_perturbed_powers} give
\[
 \|(A_\star+B_\star K)^j\|
 \leq\tau_\star((1+\gamma_\star)/2)^j,\qquad
 \|\Sigma_K\|\leq
 \frac{\sigma^2\tau_\star^2}{1-((1+\gamma_\star)/2)^2}.
\]
Thus the true closed loop is stable, and
\eqref{eq:a1_performance_difference_stationary} yields
$J(K)-J_\star\leq\|\Sigma_K\|\|S_\star\|\|K-K_\star\|_F^2
\leq C_J\varepsilon^2$.

For strong stability, use the Lyapunov argument of
\citet[Lemmas~41--42]{casselLogarithmicRegretLearning2020}
with $H=P_K^{-1/2}$, where $P_K$ is given by \eqref{eq:a1_policy_evaluation} with $G=K$ and $(A,B)=(A_\star,B_\star)$, and
$L=P_K^{1/2}(A_\star+B_\star K)P_K^{-1/2}$. It gives
\[
 \|H\|\|H^{-1}\|\leq\sqrt{\frac{J(K)}{\alpha_0\sigma^2}},
 \qquad \|L\|\leq1-\frac{\alpha_0\sigma^2}{2J(K)}.
\]
These estimates use only $Q\succeq\alpha_0I$ and $R\succeq0$.
The cited lemmas then use $R\succ0$ to bound the gain itself;
here $\|K\|\leq\bar k_\star+C_K\varepsilon_0$ supplies that bound.
\end{proof}

\section{Controller Identifiability}\label{sec:ci_appendix}\label{appendix:mvc_proof}
Recall $D_\star$ from \Cref{sec:prior_information};
 we also use $M_\star h=\Delta B_h^\top P_\star L_\star$ and write $P(\theta)$ 
 for the stabilizing Riccati solution of $\Theta(\theta)$ wherever it exists.

We use the equivalent tests below to verify controller identifiability. The gain-based CI criterion is \citet[Definition~1]{liu2026controller};
the differential reduction below is the semidefinite-cost version of
\citet[Lemma~2, Eqs.~(35)--(37)]{liu2026controller}.
\begin{lemma}
\label{lemma:controller_identifiability_tests}
Under \Cref{assumption:riccati} and \eqref{eq:a1_affine_model},
\[
 \ker\mathcal I_\star=\ker D_\star,\qquad
 \ker\mathcal H_\star=\ker\mathrm dK(\theta_\star),
\]
where $\mathrm dK(\theta_\star)$ is the differential of $\theta\mapsto K(\theta)$ at $\theta_\star$.
Consequently, CI is equivalent to
$\ker D_\star\subseteq\ker M_\star$.
\end{lemma}
\begin{proof}
Put $\Sigma_\star=\sigma^2\sum_{j\geq0}
L_\star^j(L_\star^j)^\top\succ0$.
The conditional Gaussian score gives
\[
 h^\top\mathcal I_\star h
 =\sigma^{-2}\|D_\star h\,\Sigma_\star^{1/2}\|_F^2.
\]
The smoothness of $P(\theta)$ and $K(\theta)$ near $\theta_\star$ from 
\Cref{cor:general_uniform_sensitivity} and true-plant stability
from \Cref{lemma:a1_ce} allow differentiation of the stationary
performance identity \eqref{eq:a1_performance_difference_stationary}.
Since $K(\theta_\star)=K_\star$, it gives
\[
 h^\top\mathcal H_\star h
 =2\|S_\star^{1/2}\,\mathrm dK(\theta_\star)[h]\,
                    \Sigma_\star^{1/2}\|_F^2.
\]
Here the stationary covariance depends smoothly on the stable
closed-loop matrix through its invertible Lyapunov equation.
Positivity of $S_\star$ and $\Sigma_\star$ proves the two kernel
identities. Finally, for $h\in\ker D_\star$,
\eqref{eq:general_differentiated_DARE} and invertibility of the
Lyapunov operator give $\mathrm dP(\theta_\star)[h]=0$.
Then \eqref{eq:general_differentiated_stationarity} reduces to
$S_\star\mathrm dK(\theta_\star)[h]=-M_\star h$,
and, since $S_\star\succ0$, $\mathrm dK(\theta_\star)[h]=0$ if and only if $M_\star h=0$, which proves the last assertion.
\end{proof}

\begin{proof}[Proof of \Cref{lemma:degenerate_cost_dimension}]
Under CI, every $h\in\ker D_\star$ satisfies $M_\star h=\Delta B_h^\top P_\star L_\star=0$ by \Cref{lemma:controller_identifiability_tests}, so each column of $\Delta B_h$ lies in $\ker(L_\star^\top P_\star)$, which has dimension $n-\operatorname{rank}L_\star$. Moreover, $h\mapsto\Delta B_h$ is injective on $\ker D_\star$, since $\Delta B_h=0$ forces $\Delta A_h=D_\star h=0$, hence $h=0$ by nonredundancy. Thus $\dim\ker\mathcal I_\star=\dim\ker D_\star\leq m(n-\operatorname{rank}L_\star)$, as in the proof of \citet[Corollary~2]{liu2026controller}. Then, put $Z=I-B_\star S_\star^{-1}B_\star^\top P_\star$, so that  $L_\star=ZA_\star$.
We claim $\ker Z=B_\star\ker R$. Indeed,
$ZB_\star v=B_\star S_\star^{-1}Rv=0$ if $Rv=0$.
Conversely, if $Zx=0$, let
$v=S_\star^{-1}B_\star^\top P_\star x$.
Then $x=B_\star v$ and
$S_\star v=B_\star^\top P_\star B_\star v$, hence $Rv=0$.
Moreover, $B_\star$ is injective on $\ker R$, since
$\ker R\cap\ker B_\star=\{0\}$.
The rank inequality
$\operatorname{rank}(ZA_\star)\geq
\operatorname{rank}A_\star-\dim\ker Z$
therefore proves the claimed bound.
\end{proof}

The following lemma is a quantitative adaptation of
\citet[Lemma~3 and Theorem~2]{liu2026controller}, which we extend to semidefinite input costs. We follow the same
invariance and orthogonal-decomposition argument, but use the explicit
sensitivity bounds above and then pass from optimal-feedback prediction
error to prediction error under the estimated feedback.
\begin{lemma}
\label{lemma:a1_prediction_gain}
Use the constants of \Cref{lemma:a1_ce} and
\Cref{cor:general_uniform_sensitivity}.  Assume
\begin{equation}\label{eq:a1_benign_kernel}
 \ker D_\star\subseteq\ker M_\star.
\end{equation}

If $\operatorname{rank}D_\star>0$, let
$\mu_\star=\sigma_{\min}^{+}(D_\star)>0$, the smallest nonzero
singular value, and put
\[
 \rho_{\mathrm{loc}}=\min\{r_{\mathrm{CE}},\mu_\star/(2L_{K,1})\},
 \qquad C_{\mathrm{pred}}=2L_{K,1}/\mu_\star.
\]
For $\hat\theta\in\mathbb R^p$ with $e=\hat\theta-\theta_\star$ and $\|e\|_2\leq\rho_{\mathrm{loc}}$, write $[\hat A\ \hat B]=\Theta(\hat\theta)$, $\Delta A=\hat A-A_\star$, $\Delta B=\hat B-B_\star$, and $\hat K=K(\hat\theta)$.
Then
\begin{align}
 \|\hat K-K_\star\|_F
 &\leq C_{\mathrm{pred}}\|\Delta A+\Delta B\hat K\|_F.
 \label{eq:a1_played_prediction_gain}
\end{align}
For every $h\in\ker D_\star$ with $\|h\|_2\leq r$, one has
$P(\theta_\star+h)=P_\star$ and $K(\theta_\star+h)=K_\star$.
If $D_\star=0$, then $P(\theta)=P_\star$ and $K(\theta)=K_\star$ whenever $\|\theta-\theta_\star\|_2\leq r$, and \eqref{eq:a1_played_prediction_gain} holds with
$C_{\mathrm{pred}}=0$ and $\rho_{\mathrm{loc}}=r_{\mathrm{CE}}$, without division by a singular value.
\end{lemma}
\begin{proof}
For $h\in\ker D_\star$ with $\|h\|_2\leq r$, the isometry
$[\Delta A_h\ \Delta B_h]=\mathcal Eh$ puts the perturbed model
in $\mathcal N_r$. At $P_\star$, its reference residuals satisfy
$E=D_\star h=0$ and $M=M_\star h=0$ by
\eqref{eq:a1_benign_kernel}. Hence
\eqref{eq:general_reference_gain}--\eqref{eq:general_reference_residual_identity}
give the same Riccati solution and gain, with $S_0\succ0$ by
\eqref{eq:general_reference_curvature}.
\Cref{thm:riccati_solution} identifies them uniquely.
This also proves the assertion for $D_\star=0$.

Otherwise, decompose $e=h+v$ orthogonally with
$h\in\ker D_\star$ and $v\in(\ker D_\star)^\perp$.
Since $\|h\|_2\leq\|e\|_2\leq r$, both $\theta_\star+h$ and $\hat\theta$ correspond to models in $\mathcal N_r$, so \eqref{eq:general_two_model_bounds} and exact invariance at $h$ (i.e., $K(\theta_\star+h)=K_\star$) give
\[
 \|\hat K-K_\star\|_F
 \leq L_{K,1}\|v\|_2
 \leq\frac{L_{K,1}}{\mu_\star}\|D_\star e\|_F,
\]
because $D_\star e=D_\star v$ and
$\|D_\star v\|_F\geq\mu_\star\|v\|_2$ on
$(\ker D_\star)^\perp$.
Finally,
$D_\star e=\Delta A+\Delta B\hat K-\Delta B(\hat K-K_\star)$
implies
\[
 \|\hat K-K_\star\|_F\leq\frac{L_{K,1}}{\mu_\star}
 \bigl(\|\Delta A+\Delta B\hat K\|_F
             +\|\Delta B\|\|\hat K-K_\star\|_F\bigr).
\]
Since $\|\Delta B\|\leq\|e\|_2\leq\mu_\star/(2L_{K,1})$,
the last gain term is absorbed with a factor at most $2$,
proving \eqref{eq:a1_played_prediction_gain}.
\end{proof}

\begin{proof}[Proof of \Cref{cor:mvc_known_range}]
Since $A$ is unrestricted and independent of $B$, the direction space contains all $[\Delta A\ 0]$, 
on which $D_\star h=\Delta A_h$. Thus $D_\star D_\star^*\succeq I$, hence $\sigma_{\min}^{+}(D_\star)\geq1$, 
where $D_\star^*$ is the adjoint for the Euclidean and Frobenius inner products.

For CI, every admissible $\Delta B_h$ is paired with the invisible direction $\Delta A_h=-\Delta B_hK_\star$, 
so $\ker D_\star=\{h:\Delta A_h=-\Delta B_hK_\star\}$ contains a direction for each admissible $\Delta B_h$.
 By \Cref{lemma:controller_identifiability_tests}, CI is therefore equivalent to
  $M_\star h=\Delta B_h^\top P_\star L_\star=0$, i.e., 
  $\operatorname{range}(\Delta B_h)\subseteq\ker(L_\star^\top P_\star)$, 
  for every admissible $\Delta B_h$. For $R=0$, stationarity gives $B_\star^\top P_\star L_\star=0$, 
  so $\operatorname{range}(B_\star)\subseteq\ker(L_\star^\top P_\star)$; under the range condition, $\operatorname{range}(\Delta B_h)\subseteq\operatorname{range}(B_\star)$, which proves sufficiency. 
  For necessity, recall from the proof of \Cref{lemma:degenerate_cost_dimension} that $L_\star=ZA_\star$,
   where $Z=I-B_\star S_\star^{-1}B_\star^\top P_\star$ satisfies $\ker Z=\operatorname{range}(B_\star)$ 
   when $R=0$. Hence $\operatorname{rank}L_\star=\operatorname{rank}A_\star-\dim(\operatorname{range}A_\star\cap\operatorname{range}B_\star)=\operatorname{rank}[A_\star\ B_\star]-m$.
    If $\operatorname{rank}[A_\star\ B_\star]=n$, then $\ker(L_\star^\top P_\star)$ has dimension $m$ and equals $\operatorname{range}(B_\star)$ (since $B_\star$ has full column rank),
     so CI requires $\operatorname{range}(B-B_\star)\subseteq\operatorname{range}(B_\star)$, i.e.,
      $\operatorname{range}(B)\subseteq\operatorname{range}(B_\star)$, for every $B\in\mathcal B$.

\end{proof}
\begin{remark}[Constants in the class of \Cref{cor:mvc_known_range}]\label{remark:mvc_constants}
Besides $\sigma_{\min}^{+}(D_\star)\geq1$, shown in the proof above, every system in this class satisfies $\lambda_{\min}(S_\star)\geq\alpha_0\sigma_{\min}(B_\star)^2$, since $R=0$ and $P_\star\succeq Q\succeq\alpha_0I$. Moreover, the Lyapunov argument in the proof of \Cref{lemma:a1_strong_ce} gives $\|A_\star+B_\star K_0\|\leq\sqrt{\nu_0/(\alpha_0\sigma^2)}$, so $\bar\kappa_0=\|K_0\|$ in \Cref{remark:constants} satisfies $\bar\kappa_0\leq(\sqrt{\nu_0/(\alpha_0\sigma^2)}+\vartheta)/\sigma_{\min}(B_\star)$. Hence, with $\mu_-=1$ and $s_-=\alpha_0\sigma_{\min}(B_\star)^2$, every constant in \Cref{remark:constants} and in the proof of \Cref{thm:mvc_logregret} is bounded in terms of the norm and cost bounds of \Cref{sec:prior_information} and a lower bound on $\sigma_{\min}(B_\star)$, which gives the uniform bound on $C$ stated after \Cref{cor:mvc_known_range}.
\end{remark}

\section{Proof of the Regret Upper Bound} \label{appendix:proof_mvc_regret_first_term}
\begin{remark}[Constants for \Cref{alg:mvc}]\label{remark:constants}
Use the supplied norm and cost bounds in \Cref{sec:prior_information},
and put $d=n+m$, $q=\max\{n,m\}$.
Let $\bar\kappa_0$ be a known upper bound on $\|K_0\|$ and define
\begin{equation}\label{def:kappas}
 \kappa_0=\max\!\left\{\sqrt{\frac{\nu_0}{\alpha_0\sigma^2}},
                         \bar\kappa_0,1\right\},
 \qquad \gamma_0=\frac1{2\kappa_0^2}.
\end{equation}
By the Lyapunov argument in the proof of \Cref{lemma:a1_strong_ce}, $K_0$ is
$(\kappa_0,\gamma_0)$-strongly stable. Let $s_-\leq\lambda_{\min}(S_\star)$ and $\mu_-\leq\sigma_{\min}^{+}(D_\star)$ be the supplied lower bounds, and put $\bar k_\star=\vartheta\sqrt{\nu/(\sigma^2s_-)}\geq\|K_\star\|$ (by \eqref{eq:general_elimination_bounds} and $\|P_\star\|\leq J_\star/\sigma^2$). Take $r_{\mathrm{CE}},C_K,C_J$ from \Cref{lemma:a1_ce} and $\rho_{\mathrm{loc}},C_{\mathrm{pred}}$ from \Cref{lemma:a1_prediction_gain}, evaluated with the admissible pair $(\tau_\star,\gamma_\star)$ following \Cref{lemma:general_lyapunov_inverse} and with $\nu/\sigma^2$, $\alpha_0$, $\bar k_\star$, $\tau_\star$, $\nu\tau_\star/\sigma^2$, $\vartheta$, $\alpha_1+\vartheta^2\nu/\sigma^2$, $s_-$, $\mu_-$ in place of $\|P_\star\|$, $\lambda_{\min}(P_\star)$, $\|K_\star\|$, $\|L_\star\|$, $\|P_\star L_\star\|$, $\|A_\star\|$ and $\|B_\star\|$, $\|S_\star\|$, $s_\star$, $\mu_\star$. This is valid because the restrictions in \eqref{eq:general_closed_form_radius} are upper bounds on the radius, so \Cref{lemma:general_riccati_sensitivity} holds for every smaller radius, and all constants in these lemmas are monotone in the replaced quantities. Then fix
\begin{equation}\label{def:epsilon0}
 0<\varepsilon_0\leq\min\{\rho_{\mathrm{loc}},1\},
\end{equation}
and set
\begin{equation}\label{def:kappa}
 \kappa=\max\!\left\{
  \sqrt{\frac{\nu+C_J\varepsilon_0^2}{\alpha_0\sigma^2}},
  \bar k_\star+C_K\varepsilon_0,1\right\},
 \qquad \gamma=\frac1{2\kappa^2}.
\end{equation}

For a given horizon $T\geq1$, define
\begin{equation}\label{def:xb}
\begin{gathered}
 x_b=135q\sigma^2\kappa^2\max\{\kappa_0^6(1+\vartheta)^2,4\kappa^6\}\log(5T),\\
 \lambda=\bigl(1+\max\{(1+\kappa_0)^2,\kappa^2\}\bigr)x_b,
\end{gathered}
\end{equation}
and the integer epoch length
\begin{equation}\label{def:tau}
 \tau_0=\left\lceil\frac{712(1+\kappa_0^2)\lambda d(1+\vartheta^2)}{\sigma^2\varepsilon_0^2}\right\rceil.
\end{equation}
Since $\kappa,\kappa_0\geq1$ and $\varepsilon_0\leq1$, these choices imply $\lambda\geq\max\{8n\sigma^2\log(5T),\sigma^2/20\}$ and $\tau_0\geq\gamma^{-1}\log(2\kappa)+1$; they also give $b\leq17\sigma^2\kappa^4$ for the constant $b$ in \eqref{def:a3_covariance_constants}, and hence
\begin{equation}\label{eq:tau0_bounds}
\begin{aligned}
 \tau_0-2&\geq\tfrac{1000}{9}\bigl(\log(5T^2)+d\log(2223(1+\kappa_0^2)\lambda/\sigma^2)\bigr),\\
 \tau_0-1&\geq\tfrac{1000}{9}\bigl(\log(10T^3)+n\log(2223\,b/\sigma^2)\bigr).
\end{aligned}
\end{equation}
Except in the last step of the proof of \Cref{thm:mvc_logregret}, we consider horizons satisfying $T\geq\tau_0$.
\end{remark}

\begin{lemma}
\label{lemma:a1_block_covariance}
Condition on the history at time $s$, including $x_s$ and a chosen
$(\kappa,\gamma)$-strongly stable gain $K$ on the true plant.
For an integer $N\geq1$, apply this fixed gain at all times
$s,\ldots,s+N-1$, without additional input noise, and put
\[
 X=\sum_{j=0}^{N-1}x_{s+j}x_{s+j}^\top,\qquad
 v=n\log9+\log(2/\delta),\qquad 0<\delta<1.
\]
Then, with conditional probability at least $1-\delta/2$,
\begin{equation}\label{eq:a1_block_covariance}
 X\preceq bNI,\qquad
 b=\frac{4\kappa^2\|x_s\|^2}{\gamma N}
   +\frac{8\sigma^2\kappa^2}{\gamma}
   +\frac{12\sigma^2\kappa^2v}{\gamma^2N},
\end{equation}
and, if moreover $N-1\geq\frac{1000}{9}\bigl(\log(2/\delta)+n\log(2223\,b/\sigma^2)\bigr)$, then with conditional probability at least $1-\delta$ also $X\succeq aNI$ with $a=9\sigma^2/3200$.
\end{lemma}
\begin{proof}
All probabilities below are conditional on the block-start history.
Write $F=A_\star+B_\star K$, so
$\|F^j\|\leq\kappa(1-\gamma)^j$.
For the upper bound, fix a unit vector $u$ and write the stacked
observations as $(u^\top x_{s+j})_{j=0}^{N-1}=d_u+Y_u$, where
\[
 (d_u)_j=u^\top F^jx_s,\qquad
 (Y_u)_j=\sum_{\ell=0}^{j-1}u^\top F^{j-1-\ell}w_{s+\ell}.
\]
The empty sum is zero, and $Y_u$ is centered Gaussian.
With $C_u=\operatorname{Cov}(Y_u)$, geometric sums give
\[
 \|d_u\|^2\leq\frac{\kappa^2\|x_s\|^2}{\gamma},\qquad
 \operatorname{Tr}C_u\leq\frac{\sigma^2\kappa^2N}{\gamma},\qquad
 \|C_u\|\leq\frac{\sigma^2\kappa^2}{\gamma^2}.
\]
The last bound follows by writing the noise-to-$Y_u$ map as
a sum of lag operators, which has norm at most
$\sum_{j\geq0}\|F^j\|\leq\kappa/\gamma$.
For centered Gaussian $Y$ with covariance $C$, apply
\citet[Proposition~1.1]{hsu2012tail} to $C^{1/2}g$, where $g$
is standard Gaussian. Using
$2\sqrt{\operatorname{Tr}(C^2)v}\leq
\operatorname{Tr}C+\|C\|v$ gives
\begin{equation}\label{eq:a1_gaussian_energy}
 \mathbb P\{\|Y\|^2>2\operatorname{Tr}C+3\|C\|v\}\leq e^{-v}.
\end{equation}

By \citet[Lemmas~5.2 and~5.4]{vershynin2012nonasymptotic},
there
is a deterministic set $\mathcal U$ of at most $9^n$ unit vectors such that 
$\|X\|\leq2\max_{u\in\mathcal U}u^\top Xu$.
Apply \eqref{eq:a1_gaussian_energy} to $Y_u$ for each $u\in\mathcal U$ and use
$\|d_u+Y_u\|^2\leq2\|d_u\|^2+2\|Y_u\|^2$.
A union bound gives, with probability at least $1-9^ne^{-v}=1-\delta/2$,
\[
 \|X\|\leq
 \frac{4\kappa^2\|x_s\|^2}{\gamma}
 +\frac{8\sigma^2\kappa^2N}{\gamma}
 +\frac{12\sigma^2\kappa^2v}{\gamma^2}=bN.
\]
For the lower bound, we omit $x_s$, which is fixed under the conditioning. Each of the $N-1$ states after it is conditionally Gaussian with covariance $\sigma^2I$ given its past, so its projection on any unit vector exceeds $\sigma$ in absolute value with conditional probability at least $3/10$~\citep[Eq.~(3.12)]{simchowitz2018learning}. These states thus satisfy the $(1,\Gamma_{\mathrm{sb}},3/10)$-block martingale small-ball condition of~\citet[Definition~2.1]{simchowitz2018learning} with $\Gamma_{\mathrm{sb}}=\sigma^2I$, and, as $N\leq2(N-1)$, on the event $X\preceq bNI$ their Gram matrix is at most $(N-1)\bar\Gamma$ with $\bar\Gamma=2bI\succeq\Gamma_{\mathrm{sb}}$. The condition on $N$ is \citet[Eq.~(2.7)]{simchowitz2018learning} with failure probability $\delta/2$, so the argument of~\citet[Appendix~D.1]{simchowitz2018learning} shows that, with probability at least $1-\delta/2$, either $X\not\preceq bNI$ or $X\succeq9(N-1)\sigma^2I/1600\succeq aNI$; see also~\cite{ziemann2023tutorial}.
Combining the two failure probabilities proves the claim.
\end{proof}

\begin{lemma}[cf.~{\citealp[Lemma~11]{AbbasiYadkori2011ImprovedAF}}] \label{lemma:episodic_abbasi11}
Let $\lambda > 0$ and $z_1,\ldots,z_T\in\mathbb R^d$.
Let $\tau_0<\tau_1<\cdots<\tau_N=T+1$, with epoch $i$ equal to
$\{\tau_i,\ldots,\tau_{i+1}-1\}$ for $0\leq i<N$.

Define $V_t=\lambda I_d+\sum_{s=1}^{t-1}z_sz_s^\top$ and
$Z_i=\sum_{t=\tau_i}^{\tau_{i+1}-1}z_tz_t^\top$, so that
$V_{\tau_{i+1}}=V_{\tau_i}+Z_i$.

If $c>0$ and $\lambda_{\max}(V_{\tau_i}^{-1/2} Z_i
V_{\tau_i}^{-1/2}) \le c$ for every $0\leq i<N$, then
$$\sum_{i=0}^{N-1} \sum_{t=\tau_i}^{\tau_{i+1}-1}
\|z_t\|^2_{V_{\tau_i}^{-1}}
\le \frac{c}{\log(1+c)}
\log \frac{\det(V_{T+1})}{\det(V_{\tau_0})}.$$
\end{lemma}
\begin{proof}
As in \citet[Lemma~11]{AbbasiYadkori2011ImprovedAF}, put $M_i=V_{\tau_i}^{-1/2}Z_iV_{\tau_i}^{-1/2}$, so that
$0\preceq M_i\preceq cI$. Concavity of $\log(1+x)$ gives
$x\leq c\log(1+x)/\log(1+c)$ for $0\leq x\leq c$.
Applying this inequality to the eigenvalues of $M_i$ yields
\[
\begin{aligned}
 \sum_{t=\tau_i}^{\tau_{i+1}-1}\|z_t\|_{V_{\tau_i}^{-1}}^2
 &=\operatorname{Tr}(M_i)\\
 &\leq\frac{c}{\log(1+c)}\log\det(I+M_i)\\
 &=\frac{c}{\log(1+c)}
   \log\frac{\det(V_{\tau_{i+1}})}{\det(V_{\tau_i})}.
\end{aligned}
\]
Summing over $i$ telescopes the log-determinant ratios.
\end{proof}

From now on, use the assumptions and algorithm of
\Cref{thm:mvc_logregret}, with the parameters of
\Cref{remark:constants}. Put
\[
 V_t=\lambda I_d+\sum_{s=1}^{t-1}z_sz_s^\top,
 \qquad \|M\|_{V_t}^2=\operatorname{Tr}(MV_tM^\top).
\]
Before abort, hold $\hat\Theta_t,\hat P_t,K_t$ equal to their
values at $\tau_i$ for $\tau_i\leq t<\tau_{i+1}$.
Define
\[
 D_t=\|u_t-K_\star x_t\|_{S_\star}^2,\qquad
 \Delta_{\tau_i}=\hat\Theta_{\tau_i}-\Theta_\star,\qquad
 b_0=\frac{9\varepsilon_0^2\tau_0\sigma^2}
             {3200(1+\kappa_0^2)}.
\]
For the analysis, define the constrained estimator at every epoch
boundary from the observed data, even on paths on which the algorithm has
already aborted and does not actually compute it. This latent definition
does not affect the actions and makes the following events globally defined.
We use the events below, where $\eta_t$, $t\geq\tau_0$, are auxiliary independent $\mathcal N(0,\sigma^2I_m)$ variables not used by the algorithm.
\begin{align}
\mathcal{E}_{\Theta}
&=
\left\{\begin{gathered}
\operatorname{Tr}\!\left(\Delta_{\tau_i}V_{\tau_i}\Delta_{\tau_i}^\top\right)\\
\le 4\sigma^2 n\log\!\left({5nT^2}\frac{\det(V_{\tau_i})}{\det(V_1)}\right)
+2\lambda \|\Theta_\star\|_F^2,\\
0\leq i<N_T
\end{gathered}\right\}, \label{eq:upperbound_Vt}
\\[6pt]
\mathcal E_z^{\mathrm{init}}
&=\left\{\sum_{t=2}^{\tau_0-1}z_tz_t^\top
\succeq\frac{9(\tau_0-2)\sigma^2}{3200(1+\kappa_0^2)}I_d\right\}, \label{eq:lowerbound_Vt}
\\[6pt]
\mathcal{E}_{w}
&=
\left\{
\max_{1\le t\le T}
\|w_t\|
\le
\sigma \sqrt{15 q \log(5T)}
\right\}, \label{eq:bound_process_noise}
\\[6pt]
\mathcal{E}_{\eta}
&=
\left\{
\max_{1\le t<T}
\|\eta_t\|
\le
\sigma \sqrt{15 q \log(5T)}
\right\}. \label{eq:bound_input_noise}
\end{align}

Then we have the following lemma.

\begin{lemma} \label{lemma:a3_initial_event}
Let
\[
\mathcal E_{\mathrm{init}}
=\mathcal E_\Theta\cap\mathcal E_z^{\mathrm{init}}
 \cap\mathcal E_w\cap\mathcal E_\eta.
\]
Then, under the choices in \Cref{remark:constants},
\[
\mathbb P(\mathcal E_{\mathrm{init}}^c)\leq4/(5T^2).
\]
Moreover, on $\mathcal E_w\cap\mathcal E_\eta$, $\|x_t\|^2\leq x_b$ for $t\leq\tau_0$ and $\|z_t\|^2\leq\lambda$ for $t<\tau_0$.
\end{lemma}

\begin{proof}
Let $\widetilde\Theta_t$ be the unconstrained ridge estimate. By \citet[Lemma~6]{casselLogarithmicRegretLearning2020} with $\delta=1/(5T^2)$, the inequality in \eqref{eq:upperbound_Vt} holds for $\widetilde\Theta_t$ at all $t$ with probability at least $1-1/(5T^2)$. Completing the square in \eqref{eq:constrained_ridge} shows that $\hat\Theta_t$ is the $V_t$-orthogonal projection of $\widetilde\Theta_t$ onto the affine set $\mathcal C\ni\Theta_\star$, so $\|\hat\Theta_t-\Theta_\star\|_{V_t}\leq\|\widetilde\Theta_t-\Theta_\star\|_{V_t}$, and $\mathbb P(\mathcal E_\Theta)\geq1-1/(5T^2)$.

By \citet[Lemma~34]{casselLogarithmicRegretLearning2020} with $\delta=1/(5T^2)$, $n,m\leq q$, and $\log(5T^3)\leq3\log(5T)$, $\mathbb P(\mathcal E_w),\mathbb P(\mathcal E_\eta)\geq1-1/(5T^2)$. On $\mathcal E_w\cap\mathcal E_\eta$, the states for $t\leq\tau_0$ are generated by the $(\kappa_0,\gamma_0)$-strongly stable gain $K_0$ under the disturbances $w_t+B_\star\eta_t$, whose norms are at most $\sigma(1+\vartheta)\sqrt{15q\log(5T)}$. Hence \citet[Lemma~38]{casselLogarithmicRegretLearning2020} and $\gamma_0^{-1}=2\kappa_0^2$ give $\|x_t\|^2\leq60q\sigma^2\kappa_0^6(1+\vartheta)^2\log(5T)\leq x_b$. Since $u_t=K_0x_t+\eta_t$ and $\|\eta_t\|\leq\sigma\sqrt{15q\log(5T)}\leq\sqrt{x_b}/3$, also $\|z_t\|^2\leq(1+(1+\kappa_0)^2)x_b\leq\lambda$ for $t<\tau_0$, which proves the second claim.

It remains to bound $\mathbb P(\mathcal E_w\cap\mathcal E_\eta\cap(\mathcal E_z^{\mathrm{init}})^c)$. During the warm-up, $z_t=U(x_t^\top,\eta_t^\top)^\top$ with $U=[I\ 0;\,K_0\ I]$. For $t\geq2$, given $\mathcal F_{t-1}=\sigma(x_1,u_1,\dots,x_{t-1},u_{t-1})$, the vector $(x_t,\eta_t)$ is Gaussian with covariance $\sigma^2I$; $z_1$ is excluded because $x_1=0$. 
Since $\|U^{-1}\|\leq1+\kappa_0$, $z_t$ is conditionally Gaussian with covariance $\sigma^2UU^\top\succeq\sigma^2I/(1+\kappa_0)^2$, so by~\citet[Eq.~(3.12)]{simchowitz2018learning} its projection on any unit vector exceeds $\sigma/(1+\kappa_0)$ in absolute value with conditional probability at least $3/10$. As $(1+\kappa_0)^2\leq2(1+\kappa_0^2)$, $(z_t)_{2\leq t<\tau_0}$ thus satisfies
the $(1,\Gamma_{\mathrm{sb}},3/10)$-block martingale small-ball condition of~\citet[Definition~2.1]{simchowitz2018learning} with $\Gamma_{\mathrm{sb}}=\sigma^2I/(2(1+\kappa_0^2))$, and, on $\mathcal E_w\cap\mathcal E_\eta$, its Gram matrix is at most $(\tau_0-2)\bar\Gamma$ with $\bar\Gamma=\lambda I\succeq\Gamma_{\mathrm{sb}}$. The first inequality in \eqref{eq:tau0_bounds} is \citet[Eq.~(2.7)]{simchowitz2018learning} with failure probability $1/(5T^2)$, so the argument of~\citet[Appendix~D.1]{simchowitz2018learning} gives $\mathbb P(\mathcal E_w\cap\mathcal E_\eta\cap(\mathcal E_z^{\mathrm{init}})^c)\leq1/(5T^2)$. A union bound concludes.
\end{proof}

\begin{lemma}[cf.~{\citealp[Lemma~7]{casselLogarithmicRegretLearning2020}}]\label{lemma:cassel7}
On $\mathcal E_{\mathrm{init}}$, the algorithm never aborts, and for all $1\leq t\leq T$ and $0\leq i<N_T$:
\begin{enumerate}
\item $\|x_t\|^2\leq x_b$ and $\|z_t\|^2\leq\lambda$;
\item $\|\Delta_{\tau_i}\|_{V_{\tau_i}}^2\leq b_0$ and $\|\Delta_{\tau_i}\|_F\leq\varepsilon_0$, where $V_{\tau_i}\succeq9\tau_0\sigma^2I/(3200(1+\kappa_0^2))$;
\item $K_{\tau_i}$ passes the checks of \Cref{alg:mvc} and is $(\kappa,\gamma)$-strongly stable.
\end{enumerate}
\end{lemma}
\begin{proof}
The argument follows \citet[Lemmas~7, 19, and~22]{casselLogarithmicRegretLearning2020}.
\emph{Warm-up.} By \Cref{lemma:a3_initial_event}, $\|x_t\|^2\leq x_b$ for $t\leq\tau_0$ and $\|z_t\|^2\leq\lambda$ for $t<\tau_0$. On $\mathcal E_z^{\mathrm{init}}$ and since $\lambda\geq\sigma^2/20$, $V_{\tau_i}\succeq V_{\tau_0}\succeq\frac{9\tau_0\sigma^2}{3200(1+\kappa_0^2)}I$ for all $i$.

\emph{Estimation.} Fix $i$ and suppose that $\|z_t\|^2\leq\lambda$ for all $t<\tau_i$; for $i=0$, this holds by the warm-up step. Then \citet[Lemma~37]{casselLogarithmicRegretLearning2020} gives $\log(\det V_{\tau_i}/\det V_1)\leq d\log T$, so $\mathcal E_\Theta$, $\|\Theta_\star\|_F^2\leq d\vartheta^2$, $n\leq T$, and $\lambda\geq8n\sigma^2\log(5T)$ give
\[
\begin{aligned}
 \|\Delta_{\tau_i}\|_{V_{\tau_i}}^2
 &\leq4\sigma^2n\bigl(\log(5nT^2)+d\log T\bigr)+2\lambda d\vartheta^2\\
 &\leq2\lambda d(1+\vartheta^2)\leq b_0,
\end{aligned}
\]
where the last step is \eqref{def:tau}. With the lower bound on $V_{\tau_i}$, this gives $\|\Delta_{\tau_i}\|_F\leq\varepsilon_0$. Hence, by \Cref{lemma:general_riccati_sensitivity,lemma:a1_strong_ce}, the estimated DARE has a positive definite stabilizing solution with $\hat S_{\tau_i}\succ0$, and $K_{\tau_i}$ is $(\kappa,\gamma)$-strongly stable; in particular, it passes the checks of \Cref{alg:mvc}. This proves claims 2 and 3 at epoch $i$.

\emph{No abort.} We show by induction over $i$ that $\|z_t\|^2\leq\lambda$ for all $t<\tau_i$ and that no abort occurs before $\tau_i$; the warm-up step gives the base case $i=0$. Assume this holds for $i$. By the estimation step, applied at epochs $0,\dots,i$, the state evolves from $\tau_0$ to $\tau_{i+1}$ under $(\kappa,\gamma)$-strongly stable gains that pass the checks, each held for $\tau_0\geq\gamma^{-1}\log(2\kappa)$ steps except possibly the last. Since the halving step in the proof of \citet[Lemma~39]{casselLogarithmicRegretLearning2020} is used only at controller switches, that lemma and $\mathcal E_w$ give, for $\tau_0\leq t\leq\tau_{i+1}$,
\[
 \|x_t\|\leq3\kappa\max\Bigl\{\tfrac12\|x_{\tau_0}\|,\tfrac{\kappa}{\gamma}\max_{s\leq T}\|w_s\|\Bigr\}\leq\sqrt{x_b}
\]
by \eqref{def:xb}. Hence no state-triggered abort occurs before $\tau_{i+1}$, and $\|z_t\|^2\leq(1+\kappa^2)x_b\leq\lambda$ for $\tau_0\leq t<\tau_{i+1}$, which completes the induction and proves claim 1.
\end{proof}

\begin{lemma}\label{lemma:success_event}
Let $N_i=\tau_{i+1}-\tau_i$, let
$\mathcal I_{\mathrm{full}}=\{0\leq i<N_T:N_i=\tau_0\}$, and put
\begin{equation}\label{def:a3_covariance_constants}
\begin{split}
v_T&=n\log9+\log(10T^3),\qquad a=\frac{9\sigma^2}{3200},\\
b&=\frac{4\kappa^2x_b}{\gamma\tau_0}
  +\frac{8\sigma^2\kappa^2}{\gamma}
  +\frac{12\sigma^2\kappa^2v_T}{\gamma^2\tau_0}.
\end{split}
\end{equation}
For the actual state blocks $X_i=\sum_{t=\tau_i}^{\tau_{i+1}-1}x_tx_t^\top$,
define
\begin{equation}\label{eq:a3_covariance_event}
 \mathcal E_x=\left\{\begin{gathered}
 X_i\preceq b\tau_0I\ \text{for all }0\leq i<N_T,\\
 X_i\succeq a\tau_0I\ \text{for all }i\in\mathcal I_{\mathrm{full}}
\end{gathered}\right\}.
\end{equation}
Then $\mathcal E_{\mathrm{succ}}=\mathcal E_{\mathrm{init}}\cap\mathcal E_x$ satisfies
$\mathbb P(\mathcal E_{\mathrm{succ}})\geq1-T^{-2}$.
\end{lemma}
\begin{proof}
For each $i$, let $G_i\in\mathcal F_{\tau_i}$ be the event that no abort occurred before $\tau_i$, $\|x_{\tau_i}\|^2\leq x_b$, and $K_{\tau_i}$ is $(\kappa,\gamma)$-strongly stable, and let $\widetilde X_i$ be the Gram matrix of the states obtained from $x_{\tau_i}$ by applying $K_{\tau_i}$ throughout epoch $i$, regardless of the abort rule. Conditionally on $\mathcal F_{\tau_i}$ and on $G_i$, \Cref{lemma:a1_block_covariance} with $\delta=1/(5T^3)$ shows that $\widetilde X_i$ violates the bounds in \eqref{eq:a3_covariance_event} with probability at most $1/(5T^3)$; the upper bound uses $N_i\leq\tau_0$ and $\|x_{\tau_i}\|^2\leq x_b$, and the lower bound for full epochs uses the second inequality in \eqref{eq:tau0_bounds}. On $\mathcal E_{\mathrm{init}}$, \Cref{lemma:cassel7} shows that every $G_i$ holds and $\widetilde X_i=X_i$. Hence $\mathcal E_{\mathrm{init}}\cap\mathcal E_x^c$ is contained in the union over $i$ of the events $G_i\cap\{\widetilde X_i\text{ violates the bounds}\}$, and a union bound gives $\mathbb P(\mathcal E_{\mathrm{init}}\cap\mathcal E_x^c)\leq N_T/(5T^3)\leq1/(5T^2)$. Together with \Cref{lemma:a3_initial_event}, this proves the claim.
\end{proof}

\begin{lemma}\label{lemma:a3_good_regret}
Let $a,b$ be as in \eqref{def:a3_covariance_constants}, and put
\begin{equation}\label{def:a3_potential_constant}
c=\frac{b(1+\kappa^2)(1+\kappa_0^2)}{a}.
\end{equation}
If CI holds at $\theta_\star$, then on $\mathcal E_{\mathrm{succ}}$
\begin{equation}\label{eq:a3_good_regret}
\sum_{t=\tau_0}^T D_t
\leq\|S_\star\|\frac ba C_{\mathrm{pred}}^2b_0
       \frac{c}{\log(1+c)}\,d\log(T+1)
+\mathbf1_{\{N_{N_T-1}<\tau_0\}}
       \|S_\star\|b\tau_0 C_K^2\varepsilon_0^2.
\end{equation}
If $D_\star=0$, the left side is zero.
\end{lemma}
\begin{proof}
On $\mathcal E_{\mathrm{succ}}$, \Cref{lemma:cassel7} gives no abort and
$\|\Delta_{\tau_i}\|_F\leq\varepsilon_0\leq\rho_{\mathrm{loc}}$.
Since the parametrization \eqref{eq:a1_affine_model} is isometric, this is also the parameter-error
bound required by \Cref{lemma:a1_prediction_gain}, whose assumption \eqref{eq:a1_benign_kernel} is CI by \Cref{lemma:controller_identifiability_tests}. Put
$Z_i=\sum_{t=\tau_i}^{\tau_{i+1}-1}z_tz_t^\top$ and
$\Pi_i=\Delta_{\tau_i}[I;K_{\tau_i}]$.
For each full block, the covariance bounds and that lemma give
\begin{equation}\label{eq:a1_block_loss}
\begin{split}
 \sum_{t=\tau_i}^{\tau_{i+1}-1}D_t
 &\leq b\tau_0\|S_\star\|\|K_{\tau_i}-K_\star\|_F^2\\
 &\leq b\tau_0\|S_\star\|C_{\mathrm{pred}}^2\|\Pi_i\|_F^2\\
 &\leq\|S_\star\|\frac ba C_{\mathrm{pred}}^2
             \sum_{t=\tau_i}^{\tau_{i+1}-1}\|\Delta_{\tau_i}z_t\|^2,
\end{split}
\end{equation}
where the last step uses
$\sum_t\|\Pi_ix_t\|^2=\operatorname{Tr}(\Pi_i^\top \Pi_iX_i)
\geq a\tau_0\|\Pi_i\|_F^2$ on this block.
For every epoch, the confidence bound implies
\[
 \|\Delta_{\tau_i}z_t\|^2
 \leq\|\Delta_{\tau_i}V_{\tau_i}^{1/2}\|_F^2
       \|z_t\|_{V_{\tau_i}^{-1}}^2
 \leq b_0\|z_t\|_{V_{\tau_i}^{-1}}^2.
\]
Moreover, \eqref{eq:a3_covariance_event} and \Cref{lemma:cassel7} give
\[
 Z_i=[I;K_{\tau_i}]X_i[I;K_{\tau_i}]^\top
 \preceq b\tau_0(1+\kappa^2)I_d,\qquad
 V_{\tau_i}\succeq\frac{a\tau_0}{1+\kappa_0^2}I_d.
\]
Hence $V_{\tau_i}^{-1/2}Z_iV_{\tau_i}^{-1/2}\preceq cI_d$.
Applying \Cref{lemma:episodic_abbasi11} and then
\citet[Lemma~37]{casselLogarithmicRegretLearning2020} yields
\[
\begin{aligned}
 \sum_{i=0}^{N_T-1}\sum_{t=\tau_i}^{\tau_{i+1}-1}
       \|z_t\|_{V_{\tau_i}^{-1}}^2
 &\leq\frac{c}{\log(1+c)}
       \log\frac{\det V_{T+1}}{\det V_{\tau_0}}\\
 &\leq\frac{c}{\log(1+c)}d\log(T+1).
\end{aligned}
\]
These bounds prove the full-epoch term. If the last epoch is short,
its covariance upper bound and \Cref{lemma:a1_ce} suffice:
\[
 \sum_{t=\tau_{N_T-1}}^T D_t
 \leq\|S_\star\|b\tau_0\|K_{\tau_{N_T-1}}-K_\star\|_F^2
 \leq\|S_\star\|b\tau_0C_K^2\varepsilon_0^2.
\]
Finally, if $D_\star=0$, exact invariance in
\Cref{lemma:a1_prediction_gain} makes every CE loss zero.
\end{proof}

Put
\[
 R_2=\mathbf1_{\mathcal E_{\mathrm{succ}}^c}
          \sum_{t=\tau_0}^T D_t,
 \qquad R_3=\sum_{t=1}^{\tau_0-1}D_t.
\]
\begin{lemma}\label{lemma:boundR3}
Under the assumptions of this section, there are finite constants
$C_{\mathrm{f}},C_{\mathrm{w}}$, independent of $T$, such that
\[
 \mathbb E[R_2]\leq C_{\mathrm{f}}(1+x_b/T),\qquad
 \mathbb E[R_3]\leq C_{\mathrm{w}}\tau_0.
\]
Their explicit values are given in
\eqref{def:a3_failure_constant} and \eqref{def:a3_warmup_constant}.
\end{lemma}

\begin{proof}
The argument follows \citet[Lemmas~9, 29, and~30]{casselLogarithmicRegretLearning2020}, which do not use a lower bound on $R$. Let $\tau_{\mathrm{abort}}$ be the abort time, with $\tau_{\mathrm{abort}}=T+1$ if no abort occurs, $F=\{\tau_{\mathrm{abort}}\leq T\}$, $W=\max_{s\leq T}\|w_s\|$, and $H=\max_{s<\tau_0}\|\eta_s\|$. By \Cref{lemma:cassel7}, $F\subseteq\mathcal E_{\mathrm{init}}^c$, so $\mathbb P(F)\leq T^{-2}$, and $\mathbb P(\mathcal E_{\mathrm{succ}}^c)\leq T^{-2}$ by \Cref{lemma:success_event}.

\emph{Before abort.} For $\tau_0\leq t<\tau_{\mathrm{abort}}$, the checks of \Cref{alg:mvc} give $\|K_t\|\leq\kappa$ and $\|x_t\|^2\leq x_b$, hence $D_t\leq C_{\mathrm{pre}}x_b$ with $C_{\mathrm{pre}}=\|S_\star\|(\kappa+\|K_\star\|)^2$. These terms contribute to $\mathbb E[R_2]$ at most $TC_{\mathrm{pre}}x_b\,\mathbb P(\mathcal E_{\mathrm{succ}}^c)\leq C_{\mathrm{pre}}x_b/T$.

\emph{At abort.} If $\tau_0<\tau_{\mathrm{abort}}\leq T$, the state at $\tau_{\mathrm{abort}}-1$ passed the checks, so $\|x_{\tau_{\mathrm{abort}}}\|^2\leq2\vartheta^2(1+\kappa)^2x_b+2W^2$. If $\tau_{\mathrm{abort}}=\tau_0$, \citet[Lemma~38]{casselLogarithmicRegretLearning2020} gives $\|x_{\tau_0}\|^2\leq2\kappa_0^2\gamma_0^{-2}(W^2+\vartheta^2H^2)$. Applying \citet[Lemma~35]{casselLogarithmicRegretLearning2020} to these two events, each of probability at most $T^{-2}$, and using $\log(3T^3)\leq3\log(5T)$ and $x_b\geq135q\sigma^2\log(5T)$, we get $\mathbb E[\mathbf 1_F\|x_{\tau_{\mathrm{abort}}}\|^2]\leq C_{\mathrm{a}}x_b/T^2$, where
\[
 C_{\mathrm{a}}=2\vartheta^2(1+\kappa)^2+\tfrac29\bigl(1+\kappa_0^2\gamma_0^{-2}(1+\vartheta^2)\bigr).
\]

\emph{After abort.} The fallback satisfies $\|(A_\star+B_\star K_0)^j\|\leq\kappa_0(1-\gamma_0)^j$, and the abort is decided before $w_{\tau_{\mathrm{abort}}}$ is drawn, so the subsequent noise is independent of $\{\tau_{\mathrm{abort}}=s\}$ for every $s$. Unrolling the closed loop and summing geometric series gives, with $C_y=n\sigma^2\kappa_0^2/\gamma_0$,
\[
 \mathbb E\Bigl[\mathbf1_F\sum_{t=\tau_{\mathrm{abort}}}^T\|x_t\|^2\Bigr]
 \leq\frac{2\kappa_0^2C_{\mathrm{a}}x_b}{\gamma_0T^2}+2C_yT\,\mathbb P(F),
\]
and $D_t\leq\|S_\star\|\|K_0-K_\star\|^2\|x_t\|^2$ after abort. Combining the three steps proves the bound on $\mathbb E[R_2]$ with
\begin{equation}\label{def:a3_failure_constant}
 C_{\mathrm{f}}=C_{\mathrm{pre}}+2\|S_\star\|\|K_0-K_\star\|^2
       \left(\frac{\kappa_0^2C_{\mathrm{a}}}{\gamma_0}+C_y\right).
\end{equation}

\emph{Warm-up.} The effective disturbance $w_t+B_\star\eta_t$ has covariance $\sigma^2(I+B_\star B_\star^\top)\preceq(1+\vartheta^2)\sigma^2I$. Hence, as in the proof of \citet[Lemma~29]{casselLogarithmicRegretLearning2020}, \citet[Lemma~40]{casselLogarithmicRegretLearning2020} with $x_1=0$, $Q\succeq\alpha_0I$, $R\succeq0$, and $J(K_0)\leq\nu_0$ gives $\sum_{t<\tau_0}\mathbb E\|x_t\|^2\leq(1+\vartheta^2)\nu_0\tau_0/\alpha_0$.
Since $u_t=K_0x_t+\eta_t$ and $\mathbb E\|\eta_t\|^2=m\sigma^2$, this gives $\mathbb E[R_3]\leq C_{\mathrm{w}}\tau_0$ with
\begin{equation}\label{def:a3_warmup_constant}
 C_{\mathrm{w}}=2\|S_\star\|
 \left(\|K_0-K_\star\|^2\frac{(1+\vartheta^2)\nu_0}{\alpha_0}
       +m\sigma^2\right).
\end{equation}
\end{proof}

\begin{proof}[Proof of \Cref{thm:mvc_logregret}]
The Gaussian warm-up, the state and gain bounds before abort, and the linear fallback give finite second moments at every finite time, so \eqref{eq:performance_difference_appendix} applies. Since $x_1=0$ and its terminal term is nonpositive, splitting the rounds $t\geq\tau_0$ according to $\mathcal E_{\mathrm{succ}}$ gives
\[
\mathbb E[R_T]\leq
\mathbb E\!\left[\mathbf1_{\mathcal E_{\mathrm{succ}}}\sum_{t=\tau_0}^T D_t\right]
+\mathbb E[R_2]+\mathbb E[R_3].
\]
Using \Cref{lemma:a3_good_regret,lemma:boundR3},
\begin{equation}\label{eq:a3_general_regret}
\begin{split}
\mathbb E[R_T]
&\leq\|S_\star\|\frac ba C_{\mathrm{pred}}^2b_0
       \frac{c}{\log(1+c)}\,d\log(T+1)
+\mathbf1_{\{N_{N_T-1}<\tau_0\}}
       \|S_\star\|b\tau_0 C_K^2\varepsilon_0^2\\
 &\quad+C_{\mathrm{f}}(1+x_b/T)+C_{\mathrm{w}}\tau_0,
\end{split}
\end{equation}
where $a,b,c,C_{\mathrm{f}},C_{\mathrm{w}}$ are defined in
\eqref{def:a3_covariance_constants}, \eqref{def:a3_potential_constant},
\eqref{def:a3_failure_constant}, and \eqref{def:a3_warmup_constant}.
The constants $\varepsilon_0$, $\kappa$, $\gamma$, $C_K$, $C_{\mathrm{pred}}$, $C_{\mathrm{f}}$, and $C_{\mathrm{w}}$ do not depend on $T$, whereas $x_b$ and $\lambda$ are proportional to $\log(5T)$ by \eqref{def:xb}, and so are $\tau_0$ and $b_0$ up to the rounding in \eqref{def:tau}. Since $v_T\leq(n\log9+3)\log(5T)$, the constants $b$ and $c$, which depend on $T$ only through $x_b/\tau_0$ and $v_T/\tau_0$, are bounded uniformly in $T$, and so is $c/\log(1+c)$, which is increasing in $c$. As $1\leq\log(T+1)\leq\log(5T)$, every term of \eqref{eq:a3_general_regret} is therefore at most a $T$-independent multiple of $\log^2(5T)$.

Finally, if $T<\tau_0$, only warm-up inputs are played up to time $T$. Since $D_t\geq0$, \eqref{eq:performance_difference_appendix} and the warm-up step in the proof of \Cref{lemma:boundR3}, which does not use $T\geq\tau_0$, give $\mathbb E[R_T]\leq C_{\mathrm w}\tau_0$, which is again at most a $T$-independent multiple of $\log^2(5T)$.
\end{proof}

\section{Proof of the Regret Lower Bound}
\label{appendix:mvc_lower_bound_proof}
In the family \eqref{eq:lower_affine_family}, take the two alternatives
\begin{equation}\label{eq:system_lower_bound}
 \theta=\nu\Delta,\qquad \Delta=\tfrac12T^{-1/4},\qquad
 \nu\sim\operatorname{Unif}\{-1,1\},\qquad B_\nu:=B_{\nu\Delta}.
\end{equation}
For $T\geq10^4$ both belong to the stated parameter interval.
By Yao's principle~\citep{Yao1977ProbabilisticCT}, as in~\cite{casselLogarithmicRegretLearning2020}, \Cref{thm:lower_bound} follows from the following bound for deterministic policies under the random sign $\nu$.
\begin{theorem}[cf.~{\citealp[Theorem~13]{casselLogarithmicRegretLearning2020}}]\label{thm:cassel13}
For $T\geq10^4$, every deterministic learning algorithm for the alternatives
\eqref{eq:system_lower_bound} satisfies
\[
 \mathbb E[R_T]\geq\frac{3.6\cdot10^{-4}}{1+r}\sqrt T-2.4,
\]
where the expectation includes the process noise and the random sign $\nu$.
\end{theorem}

For $P=\operatorname{diag}(1,p)$, one has $A^\top PA=\operatorname{diag}(0,1/16)$ and $B_\nu^\top PA=[0,\ \nu\Delta/4]$, so the DARE reduces to
\[
 p=\frac{17}{16}-\frac{\Delta^2}{16\beta},\qquad
 \beta=r+\Delta^2+p,
\]
whose unique positive root gives the stabilizing solution, with gains $K_\nu=[0,-\nu k_\star]$, $k_\star=\Delta/(4\beta)$, and closed loops $L_\nu=A+B_\nu K_\nu$ with eigenvalues $\{0,-\nu k_\star\}$. Since $P=Q+K_\nu^\top RK_\nu+L_\nu^\top PL_\nu\succeq I_2$ and $\Delta\leq0.05$,
\begin{equation}\label{eq:lower_model_bounds}
\begin{gathered}
 1\leq p\leq17/16,\qquad 1+r\leq\beta\leq1.065+r,\\
 \frac{0.23\Delta}{1+r}\leq k_\star\leq\frac{0.25\Delta}{1+r},\qquad
 \|L_\nu\|^2\leq0.065,
\end{gathered}
\end{equation}
and $\operatorname{Tr}(P)=1+p\leq2.0625$, $\|B_\nu\|^2=1+\Delta^2\leq1.01$.

\begin{lemma}[cf.~{\citealp[Lemma~14]{casselLogarithmicRegretLearning2020}}]\label{lemma:cassel14}
Fix $\nu\in\{-1,1\}$ and any causal policy, and let $\mathbb E_\nu$ denote the expectation conditional on this sign. Then
\[
 \sum_{t=1}^{T-1}\mathbb E_\nu[(u_t-K_\nu x_t)^2]\leq\frac{1.2\,(\mathbb E_\nu[R_T]+2.4)}{1+r}.
\]
\end{lemma}

\begin{proof}
If $\mathbb E_\nu[R_T]=\infty$, the claim is immediate. Otherwise $x_t$, $t\leq T$, and, since $B_\nu$ has full column rank, $u_t$, $t\leq T-1$, have finite second moments. Put $v_t=u_t-K_\nu x_t$. Applying \eqref{eq:performance_difference_appendix} over the first $T-1$ rounds and adding the last stage cost minus $J_\star=\operatorname{Tr}(P)$ gives
\begin{equation}\label{eq:lemma14_eq2}
 \mathbb E_\nu[R_T]
 =\beta\sum_{t=1}^{T-1}\mathbb E_\nu[v_t^2]
 -\mathbb E_\nu[x_T^\top(P-I_2)x_T]
 +\mathbb E_\nu[ru_T^2]-\operatorname{Tr}(P).
\end{equation}
Unrolling the dynamics, using $\|L_\nu^j\|_F^2\leq2\|L_\nu\|^{2j}$, and applying Cauchy--Schwarz as in
\citet[Lemma~32]{casselLogarithmicRegretLearning2020} gives
\[
\begin{aligned}
 \mathbb E_\nu\|x_T\|^2
 &\leq\frac{2\|B_\nu\|^2}{1-\|L_\nu\|^2}\sum_{t=1}^{T-1}\mathbb E_\nu[v_t^2]
       +\frac{4}{1-\|L_\nu\|^2}\\
 &\leq2.2\sum_{t=1}^{T-1}\mathbb E_\nu[v_t^2]+4.4.
\end{aligned}
\]
Since $0\preceq P-I_2\preceq0.0625I_2$, $\beta\geq1+r$, $ru_T^2\geq0$, and
$\operatorname{Tr}(P)\leq2.0625$, inserting this in \eqref{eq:lemma14_eq2} yields
\begin{equation}\label{eq:lemma14_final}
 \mathbb E_\nu[R_T]\geq0.86(1+r)\sum_{t=1}^{T-1}\mathbb E_\nu[v_t^2]-2.4,
\end{equation}
because $\beta-0.1375\geq0.86(1+r)$ and $0.275+2.0625\leq2.4$. In particular $\mathbb E_\nu[R_T]\geq-2.4$, and rearranging proves the claim since $1/0.86\leq1.2$.
\end{proof}

\begin{proof}[Proof of \Cref{thm:cassel13}]
We proceed similarly to the proof of 
\citet[Theorem~13]{casselLogarithmicRegretLearning2020}.
Infinite expected cost is immediate, so assume finite cost and,
for a contradiction, $\mathbb E[R_T]<c_r\sqrt T-2.4$, where
$c_r=3.6\cdot10^{-4}/(1+r)$. Put $v_t=u_t-K_\nu x_t$ and
let $\mathbb P_\nu$ be the trajectory law through $x_{T-1}$.
As in \citet[Lemma~15]{casselLogarithmicRegretLearning2020},
the Gaussian KL-divergence chain rule gives
$\mathrm{KL}(\mathbb P_\nu\|\mathbb P_{-\nu})
=2\Delta^2\sum_{t=1}^{T-2}\mathbb E_\nu[u_t^2]$.
Pinsker in both directions, followed by averaging over the sign, gives
\begin{equation}\label{eq:TV_bound}
 \operatorname{TV}(\mathbb P_1,\mathbb P_{-1})^2
 \leq\Delta^2\sum_{t=1}^{T-1}\mathbb E[u_t^2].
\end{equation}
Using $u_t^2\leq2v_t^2+2k_\star^2\|x_t\|^2$,
\Cref{lemma:cassel14}, and the cost bound
$\sum_{t<T}\mathbb E\|x_t\|^2\leq\mathbb E[R_T]+2.1T$ (which holds since $ru_t^2\ge 0$), we obtain
\[
\begin{aligned}
 \operatorname{TV}(\mathbb P_1,\mathbb P_{-1})^2
 &\leq2\Delta^2\Bigl(\frac{1.2(\mathbb E[R_T]+2.4)}{1+r}
 +k_\star^2(\mathbb E[R_T]+2.1T)\Bigr)\\
 &\leq\frac1{16},
\end{aligned}
\]
where the last inequality holds because, by $T\geq10^4$, 
$\mathbb E[R_T]\leq c_r\sqrt T$, and \eqref{eq:lower_model_bounds}, the first term is at most $0.6c_r/(1+r)+1.44T^{-1/2}\leq0.015$ and the second at most $(c_rT^{-1/2}+2.1)/128\leq0.017$.

For each $2\leq t<T$, put $A_t=\{u_tx_{2,t}\geq0\}$ and
$C_t=\{|x_{2,t}|>2/5\}$. These are events of the same observed
trajectory under both signs because the policy is deterministic.
Since $x_{2,t}=u_{t-1}+w_{2,t-1}$ under either sign, its conditional Gaussian density
is at most $(2\pi)^{-1/2}$, giving
$\mathbb P_{-1}(C_t)\geq1-4/(5\sqrt{2\pi})>2/3$.
On $A_t\cap C_t$ under $\nu=1$, or $A_t^c\cap C_t$ under $\nu=-1$,
the control has the wrong sign and $v_t^2=(u_t + \nu k_\star x_{2,t})^2\geq4k_\star^2/25$.
Consequently,
\[
\begin{aligned}
 \mathbb E[v_t^2]
 &\geq\frac{2k_\star^2}{25}
    \bigl(\mathbb P_1(A_t\cap C_t)+\mathbb P_{-1}(A_t^c\cap C_t)\bigr)\\
 &\geq\frac{2k_\star^2}{25}
    \bigl(\mathbb P_{-1}(C_t)-\operatorname{TV}(\mathbb P_1,\mathbb P_{-1})\bigr)
 \geq\frac{k_\star^2}{30}.
\end{aligned}
\]
Summing and averaging \eqref{eq:lemma14_final} over the two signs now gives
\[
 \mathbb E[R_T]\geq\frac{0.86(1+r)}{30}(T-2)k_\star^2-2.4
 \geq c_r\sqrt T-2.4,
\]
by using \eqref{eq:lower_model_bounds} and that $0.86\cdot0.23^2(T-2)/(120\sqrt T)\geq3.6\cdot10^{-4}\sqrt T$ for $T\geq10^4$. This is a contradiction.
\end{proof}

\paragraph{Relation to CI.} Replacing $\Delta$ by $\theta$ in the DARE solution above gives $K(\theta)=[0,-k(\theta)]$ with $k(\theta)=\theta/[4(r+\theta^2+p_\theta)]$, so $k(0)=0$ and $k(\theta)\sim\theta/[4(r+17/16)]$ as $\theta\to0$. Under a gain $[0,-k]$ with $|k|<1$, the state obeys $x_{1,t+1}=(1/4-\theta k)x_{2,t}+w_{1,t}$ and $x_{2,t+1}=-kx_{2,t}+w_{2,t}$, so $x_2$ has stationary variance $1/(1-k^2)$. Since only $\theta$ is unknown, the derivative of the mean of $x_{t+1}$ with respect to $\theta$ is $(-kx_{2,t},0)^\top$, and hence $\mathcal I_\theta(K(\theta))=k(\theta)^2/(1-k(\theta)^2)$. At $\theta=0$, these recursions give $J_{\Theta(0)}([0,-k])=1+(17/16+rk^2)/(1-k^2)=33/16+(r+17/16)k^2+O(k^4)$; with $k(0)=0$ and $k'(0)=1/(4(r+17/16))$, this yields $\mathcal H_0=2(r+17/16)k'(0)^2=1/[8(r+17/16)]>0=\mathcal I_0$. Finally, when $\theta$ is the true parameter, $D_\star h=\Delta B_hK(\theta)$ has the single nonzero entry $-hk(\theta)$, so $\sigma_{\min}^{+}(D_\star)=|k(\theta)|$.

\end{document}